\documentclass[a4paper,10pt]{article}

\usepackage{amsmath,amssymb}
\usepackage[T1]{fontenc} 
\usepackage[frenchb]{babel}

\newcommand{\Cset}{\mathbf{C}}
\newcommand{\Fset}{\mathbf{F}}

\newcommand{\Qset}{\mathbf{Q}}
\newcommand{\Rset}{\mathbf{R}}
\newcommand{\Zset}{\mathbf{Z}}

\newcommand{\car}{\lambda}
\newcommand{\cqi}{corps quadratique imaginaire}
\newcommand{\Etors}[1][p]{E\lbrack #1\rbrack}
\newcommand{\Exc}[1][E]{\mathrm{Exc}(#1)}
\newcommand{\Galbar}[1][K]{\mathrm{Gal}(\Qbar/{#1})}

\newcommand{\ideal}[1][L]{\mathfrak{q}}
\newcommand{\ide}[1][K]{I_{#1}}
\newcommand{\integers}[1][K]{\mathcal{O}_{#1}}

\newcommand{\Qbar}{\overline{\mathbf{Q}}}
\newcommand{\rec}[1]{{r}_{#1}}
\newcommand{\rep}[1][p]{\rho_{#1}}
\newcommand{\Res}{\mathrm{Res}}

\newcommand{\ssi}{si et seulement si }

\newtheorem{theorem}{Théorème}
\newtheorem{lemma}{Lemme}[section]
\newtheorem{proposition}{Proposition}[section]

\newtheorem{corollary}{Corollaire}[section]

\newtheorem{example}{Exemple}[section]

\newcommand\myprod{\mathop{\mathchoice
{\text{\huge$\ast$}} 
{\text{\Large$\ast$}} 
{\text{\normalsize$\ast$}} 
{\text{\footnotesize$\ast$}} 
}\displaylimits}

\newtheorem{notation}{Notation}

\begin{document}

\title{Crit\`eres d'irr\'eductibilit\'e pour les repr\'esentations \\ 
des courbes elliptiques}

\author{Nicolas Billerey}

\date{}

\maketitle

\begin{abstract}
Let $E$ be an elliptic curve defined over a number field $K$. We say that a prime number $p$ is
exceptional for $(E,K)$ if $E$ admits a $p$-isogeny defined over $K$. The so-called 
exceptional set of all such prime numbers is finite if and only if $E$ does not have complex multiplication 
over~$K$. In this paper, we prove that the exceptional set is 
included in the set of prime divisors of an explicit list of integers (depending on $E$ and $K$), whose infinitely many of them are non-zero. It
provides an efficient algorithm for computing it in the finite case. Other less general but rather 
useful criteria are given, as well as several numerical examples.
\end{abstract}

			\section*{Introduction}

Soient $\overline{\Qset}$ la cl\^oture alg\'ebrique de $\Qset$ dans $\Cset$ et $K$ un corps de nombres contenu dans $\Qbar$. \'Etant donn\'es une courbe elliptique $E$ d\'efinie sur $K$ et un nombre premier $p$, on note $\Etors$ le groupe des points de $p$-torsion de la courbe~$E$. C'est un espace vectoriel de dimension~$2$ sur le corps $\Fset_p=\Zset/p\Zset$ muni d'une action du groupe de Galois $\mathrm{G}_K=\Galbar$. Cela fournit un homomorphisme 
\begin{equation*}
\rep: \mathrm{G}_{K}\longrightarrow \mathrm{Aut}(\Etors)\simeq \mathrm{GL}_2(\Fset_p).
\end{equation*}	
Serre a d\'emontr\'e (\cite{Ser72}) que si $E$ est sans multiplication complexe sur $\Qbar$, il existe une constante $c(E,K)$ telle que pour tout nombre premier $p>c(E,K)$, la repr\'esentation $\rep$ est surjective. Il a \'egalement pos\'e la question (toujours ouverte, y compris pour $K=\Qset$) 
de savoir si $c(E,K)$ peut \^etre choisie ind\'ependamment de $E$ (\cite{Ser79}). 

Dans ce travail, on s'intéresse à l'ensemble, noté $\Exc[E/K]$, des nombres premiers~$p$ pour lesquels la représentation $\rep$ ci-dessus est réductible. On dit alors que $p$ est exceptionnel pour le couple $(E,K)$. L'ensemble $\Exc[E/K]$ est généralement fini. Plus précisément, $\Exc[E/K]$ est fini si et seulement $E$ n'a pas de multiplication complexe sur $K$, i.e. {$\mathrm{End}_K(E)=\Zset$} (prop.~\ref{prop:CNS_finitude}). On s'intéresse ici à la question suivante.

\medskip

\noindent{\it Question~1. }
Le corps $K$ et la courbe $E$ \'etant donn\'es, comment déterminer \emph{explicitement} l'ensemble $\Exc[E/K]$ lorsqu'il est fini?

\medskip

Lorsque $E$ est sans multiplication complexe, Pellarin (\cite{Pel01}), \`a la suite de Masser et W\"ustholz, a obtenu, comme corollaire de ses travaux, une majoration explicite des nombres premiers exceptionnels. Cependant, en raison des constantes qui y appara\^issent, ce r\'esultat ne se prête malheureusement pas à une détermination \emph{explicite} de l'ensemble des nombres premiers exceptionnels. En utilisant des arguments de théorie du corps de classes, on obtient un critère (th\'eor\`eme~\ref{th:critere2}) portant sur la réduction de $E$ en chaque place finie de $K$ offrant une r\'eponse satisfaisante \`a la question~$1$. De plus, dans le cas o\`u $K$ est de degr\'e impair sur~$\Qset$, on dispose d'un autre r\'esultat (th\'eor\`eme~\ref{th:critere1}) particuli\`erement simple \`a mettre en application.

Ces th\'eor\`emes sont illustrés numériquement dans le \S\ref{s:exemples} où l'on détermine explicitement l'ensemble des nombres premiers exceptionnels de plusieurs courbes elliptiques.

\medskip

On s'intéresse également dans ce travail à la question suivante.

\medskip

\noindent{\it Question~2. }
Soient $K$ un corps de nombres et $\mathcal{E}$ un ensemble infini de courbes elliptiques d\'efinies sur $K$ tels que pour toute courbe $E$ de l'ensemble $\mathcal{E}$, $\Exc[E/K]$ soit fini. Peut-on trouver une constante uniforme $\alpha(\mathcal{E},K)$ telle que pour toute courbe elliptique $E$ appartenant \`a $\mathcal{E}$, la repr\'esentation $\rep$ soit irr\'eductible d\`es que $p>\alpha(\mathcal{E},K)$?

\medskip

Dans le cas o\`u $\mathcal{E}$ est l'ensemble de toutes les courbes elliptiques sans multiplication complexe d\'efinies sur~$K$, cette question est une \'etape (importante) vers la r\'esolution de la question uniforme de Serre (voir \cite{BiPa08} pour plus de d\'etails et de nouvelles avanc\'ees). Lorsque $K=\Qset$ et $\mathcal{E}$ est l'ensemble de toutes les courbes elliptiques d\'efinies sur $\Qset$, Mazur a montr\'e (\cite{Maz78}) que tel est le cas avec $\alpha(\mathcal{E},\Qset)=163$. Dans le cas o\`u $\mathcal{E}$ est l'ensemble des courbes semi-stables, Kraus a obtenu des r\'esultats uniformes et effectifs pour diff\'erentes familles corps de nombres, notamment les corps quadratiques et cubiques (\cite{Kra96,Kra07}). Dans ce travail, on g\'en\'eralise aux corps de nombres plusieurs résultats connus sur~$\Qset$ (prop.~\ref{prop:Phi_1} et~\ref{prop:Phi_2}). On obtient ainsi quelques r\'eponses dans la direction de la question~$2$ pour des ensembles $\mathcal{E}$ de courbes elliptiques ayant mauvaise réduction additive en une place finie de $K$ et un \og défaut de semi-stabilité\fg\ particulier. Ces résultats sont particulièrement utiles d'un point de vue numérique et sont illustrés au \S\ref{s:exemples}.

	\section{\'Enonc\'es des r\'esultats}\label{s:enonces_resultats}

	\subsection{Une loi utile}\label{ss:Notations}

On note $M_{\Zset}$ le sous-ensem\-ble de $\Zset[X]$ constitué des polynômes unitaires ne s'annulant pas en~$0$. L'application 
\begin{equation*}
\begin{array}{ccl}
M_{\Zset}\times M_{\Zset} & \longrightarrow & \Zset[X] \\
(P,Q) & \longmapsto & (P\ast Q)(X)=\Res_Z\left(P(Z),Q\left(X/Z\right)Z^{\deg Q}\right)
\end{array}
\end{equation*}
où $\Res_Z$ désigne le résultant par rapport à la variable~$Z$, définit une loi de monoïde commutatif 
sur~$M_{\Zset}$ d'élément neutre $X-1$ (lemme~\ref{lem:loi_monoide}). De plus, les racines complexes 
de $P\ast Q$ sont exactement les produits d'une racine de $P$ par une racine de $Q$, comptées avec 
multiplicités (loc.~cit.). Concrètement, si
\begin{equation*}
P(X)=\prod_{i}(X-\alpha_i)\quad\textrm{et}\quad Q(X)=\prod_j(X-\beta_j),
\end{equation*}
alors, on a $(P\ast Q)(X)=\prod_{i,j}(X-\alpha_i\beta_j)$.

\'Etant donn\'es $P\in M_{\Zset}$ et $k\geq 1$, on convient de noter 
\begin{equation*}
P^{* k}=\underbrace{P\ast\cdots\ast P}_{k\textrm{ fois}} \quad\textrm{et}\quad P^{* 0}(X)=X-1.
\end{equation*}
Par ailleurs, soient $P\in M_{\Zset}$ et $r\geq 1$. Il existe alors un unique polynôme $P^{(r)}\in M_{\Zset}$ tel que
\begin{equation*}
P^{(r)}(X^r)=(P\ast \Psi_r)(X)
\end{equation*}
où $\Psi_r(X)=X^r-1$ (lemme~\ref{lem:morphisme_monoide}). Les racines complexes de $P^{(r)}$ sont 
exactement les puissances $r$-ièmes des racines complexes de $P$ comptées avec multiplicités. Autrement
dit, si 
\begin{equation*}
P(X)=\prod_{i}(X-\alpha_i),\quad\textrm{alors}\quad P^{(r)}(X)=\prod_{i}(X-\alpha_i^r).
\end{equation*}
Enfin, l'application $P\mapsto P^{(r)}$ est un morphisme de monoïdes pour la loi~$\ast$.

	\subsection{Résultats}

On \'enonce deux r\'esultats dans la direction de la question~$1$ et deux dans la direction de la question~$2$ que l'on illustre dans le \S\ref{s:exemples} sur des exemples concrets.

On fixe un corps de nombres $K$ contenu dans $\Qbar$ et une courbe elliptique $E$ d\'efinie sur $K$. On note $d$ le degré de $K$ sur $\Qset$, $D_K$ son discriminant, $\integers$ son anneau d'entiers, $h$ son nombre de classes et $N_{K/\Qset}$ la norme de l'extension $K/\Qset$.

Soit $\mathfrak{q}$ un id\'eal premier de $\integers$ en lequel $E$ a bonne r\'eduction. On pose 
\begin{equation*}
P_{\mathfrak{q}}(X)=X^2-t_{\mathfrak{q}}X+N(\mathfrak{q})\in\Zset[X]
\end{equation*}
o\`u $N(\mathfrak{q})$ est le cardinal du corps r\'esiduel $\integers/\mathfrak{q}$ et
\begin{equation*}
t_{\mathfrak{q}}=N(\mathfrak{q})+1-A_{\mathfrak{q}},
\end{equation*}
avec $A_{\mathfrak{q}}$ le nombre de points sur le corps $\integers/\mathfrak{q}$ de la r\'eduction de $E$ en~$\mathfrak{q}$. 

\medskip

Soit $\ell$ un nombre premier et 
\begin{equation*}
\ell\integers=\prod_{\ideal\mid\ell}\ideal^{v_{\ideal}(\ell)}
\end{equation*}
la décomposition de $\ell\integers$ en produit d'idéaux premiers de~$\integers$. On suppose que $E$ a 
bonne réduction en chaque idéal premier au-dessus de $\ell$. Par abus de langage, on dit alors que $E$ a bonne r\'eduction en~$\ell$. Dans ce cas, on associe à $\ell$ (on rappelle
que $E$ et $K$ sont fixés) un polynôme $P_{\ell}^*$ à coefficients entiers dont certaines valeurs 
spéciales vont permettre de déterminer essentiellement l'ensemble $\Exc[E/K]$. Précisément, on pose~:
\begin{equation}\label{eq:polynôme_P_ell*}
P_{\ell}^*=\myprod_{\ideal\mid\ell} \left(P_{\ideal}^{(12v_{\ideal}(\ell))}\right)\in\Zset[X],
\end{equation}
où, d'une part le produit $\ast$, pris au sens de la définition du \S\ref{ss:Notations}, porte sur tous les idéaux
premiers de $\integers[K]$ divisant~$\ell$ et, d'autre part les exposants $(12v_{\ideal}(\ell))$ renvoient 
à la notation adoptée également au \S\ précédent. Ce polynôme ne dépend que de la famille de triplets d'entiers 
$\{(t_{\ideal},v_{\ideal}(\ell),f_{\ideal})\}_{\ideal\mid\ell}$ où $N(\ideal)=\ell^{f_{\ideal}}$. 
Ses racines complexes sont de module~$\ell^{6d}$ (lemme~\ref{lem:polynome_P_ell}). On considère alors
l'entier (essentiel dans la suite)~:
\begin{equation*}
B_{\ell}=\prod_{k=0}^{\left[\frac{d}{2}\right]}P_{\ell}^*\left(\ell^{12k}\right)
\end{equation*}
où $[d/2]$ désigne la partie entière de $d/2$. 

\medskip

Dans la direction de la question~$1$ de l'\emph{Introduction}, on montre le critère suivant.
\begin{theorem}\label{th:critere1}
Soit $p$ un nombre premier exceptionnel pour $(E,K)$. Alors, on est dans l'une des situations suivantes~:
\begin{enumerate}
	\item $p$ divise $6D_K$;
	\item il existe un idéal premier $\mathfrak{p}$ de~$\integers$ divisant~$p$ en lequel $E$ a mauvaise réduction additive avec potentiellement bonne réduction supersingulière.
	\item pour tout nombre premier $\ell$ de bonne réduction, le nombre premier $p$ divise l'entier $B_{\ell}$ (si $d=1$, on suppose $\ell\not=p$).
\end{enumerate}
\end{theorem}

\medskip

Supposons que $E$ soit donnée par une équation de Weierstrass à coefficients dans l'anneau $\integers$ de discriminant $\Delta_E$. On déduit du théorème~\ref{th:critere1} le corollaire suivant.
\begin{corollary}\label{cor:discriminant}
Soit $p$ un nombre premier exceptionnel pour $(E,K)$. Alors, on est dans l'une des situations suivantes~:
\begin{enumerate}
	\item $p$ divise $6D_K N_{K/\Qset}(\Delta_E)$;
	\item pour tout nombre premier $\ell$ de bonne réduction, le nombre premier $p$ divise l'entier $B_{\ell}$ (si $d=1$, on suppose $\ell\not=p$).
\end{enumerate} 
\end{corollary}

Les racines complexes de $P_{\ell}^*$ étant de module~$\ell^{6d}$ (lemme~\ref{lem:polynome_P_ell}), on a en particulier~:
\begin{equation*}
d\textrm{ impair }\Longrightarrow B_{\ell}\not=0.
\end{equation*}
On déduit alors du théorème~\ref{th:critere1} le corollaire suivant.
\begin{corollary}[cas du degré impair]\label{cor:degre_impair}
On suppose que l'extension $K/\Qset$ est de degré $d$ impair. Alors, l'ensemble des nombres premiers exceptionnels pour $E$ est \emph{fini}. De plus, si $p$ un nombre premier exceptionnel pour $(E,K)$, alors, pour tout nombre premier $\ell$ de bonne réduction, le nombre premier $p$ divise l'entier \emph{non nul}
\begin{equation*}
6D_KN_{K/\Qset}(\Delta_E)B_{\ell}.
\end{equation*}
\end{corollary}

\medskip 

\noindent{\it Remarque. }Plus généralement, on démontre à partir de la proposition~\ref{prop:CNS_finitude} et de résultats classiques de la théorie de la multiplication complexe que les propriétés suivantes sont équivalentes~:
\begin{enumerate}
	\item le corps $K$ ne contient pas le corps de classes de Hilbert d'un \cqi;
	\item pour toute courbe elliptique $E$ définie sur $K$, l'ensemble $\Exc[E/K]$ est fini.
\end{enumerate}

\medskip 

La situation est plus compliquée dans le cas des extensions de degré pair. Bien que le crit\`ere du th\'eor\`eme~\ref{th:critere1} ci-dessus s'applique toujours (cf.~\S\ref{ss:cas_quad}), on n'a plus la garantie, pour une courbe ayant un ensemble exceptionnel fini, qu'il existe un nombre premier $\ell$ de bonne r\'eduction pour lequel l'entier $B_{\ell}$ soit non nul (comme le montre l'exemple~\ref{ex:biquad}). L'\'enonc\'e plus g\'en\'eral ci-dessous permet de contourner cette difficult\'e et constitue notre r\'esultat principal en direction de la question~$1$.
\begin{theorem}\label{th:critere2}
Soit $p$ un nombre premier exceptionnel pour $(E,K)$. Alors, on est dans l'une des situations suivantes~:
\begin{enumerate}
	\item $p$ divise $6D_K$;
	\item il existe un idéal premier $\mathfrak{p}$ de~$\integers$ divisant~$p$ en lequel $E$ a mauvaise réduction additive avec potentiellement bonne réduction supersingulière.
	\item pour tout id\'eal premier $\ideal$ de bonne réduction, le nombre premier $p$ divise l'entier 
\begin{equation*}
R_{\ideal}=\prod_{k=0}^{\left[\frac{d}{2}\right]}\Res\left(P_{\ideal}^{(12h)},\left(\mathfrak{m}_{\gamma_{\ideal}}^{(12)}\right)^{* k}\right),
\end{equation*}
o\`u $\ideal^h=\gamma_{\ideal}\mathcal{O}_K$ et $\mathfrak{m}_{\gamma_{\ideal}}$ est le polyn\^ome minimal de $\gamma_{\ideal}$ sur~$\Qset$ (si $d=1$, on suppose que $\ideal$ ne divise pas~$p$).
\end{enumerate}
De plus, si $E$ est sans multiplication complexe sur~$\overline{\Qset}$, alors $R_{\ideal}\not=0$ pour une infinit\'e d'id\'eaux premiers~$\ideal$.
\end{theorem}

\bigskip

Dans la direction de la question~$2$ de l'\emph{Introduction}, on généralise aux corps de nombres plusieurs résultats connus sur~$\Qset$. Soit $\mathfrak{q}$ un idéal premier de $\integers$ de caractéristique résiduelle~$\ell$. On a
\begin{equation*}
N(\mathfrak{q})=\lvert\integers/\mathfrak{q}\rvert=\ell^{f_{\mathfrak{q}}},
\end{equation*}
où $f_{\mathfrak{q}}$ est le degré résiduel de $\mathfrak{q}$. On suppose que $E$ a mauvaise réduction additive en~$\mathfrak{q}$ avec potentiellement bonne réduction. Alors, pour tout nombre premier $p\geq 3$ tel que $p\not=\ell$, l'action de $I_{\mathfrak{q}}$, sous-groupe d'inertie en $\mathfrak{q}$, sur $\Etors$ se fait par l'intermédiaire d'un certain quotient fini $\Phi_{\mathfrak{q}}$ de~$I_{\mathfrak{q}}$ (\cite{SeTa68})~:
\begin{equation*}
I_{\mathfrak{q}}\longrightarrow \Phi_{\mathfrak{q}}\hookrightarrow \mathrm{Aut}(\Etors).
\end{equation*}
On a alors les deux r\'esultats suivants. Ceux-ci sont connus pour $K=\Qset$ et ont \'et\'e utilis\'es par Serre dans \cite[\S5]{Ser72} pour traiter des exemples num\'eriques.
\begin{proposition}\label{prop:Phi_1}
On suppose que le groupe $\Phi_{\mathfrak{q}}$ n'est pas cyclique. Alors, la repr\'esentation $\rep$ est irr\'eductible pour tout nombre premier $p\geq 5$.
\end{proposition}

\begin{proposition}\label{prop:Phi_2}
On suppose que pour tout entier $n\geq 0$, l'ordre du groupe $\Phi_{\mathfrak{q}}$ ne divise pas $N(\mathfrak{q})^n(N(\mathfrak{q})-1)$. Alors, la repr\'esentation $\rep$ est irr\'eductible pour tout nombre premier $p\geq 3$ tel que $p\not=\ell$.
\end{proposition}
\noindent{\it Remarque. }\label{rem:precisionpaumoins5}
Si $\ell\geq 5$, on peut remplacer l'hypothèse par~: l'ordre du groupe $\Phi_{\mathfrak{q}}$ ne divise pas $N(\mathfrak{q})-1$.

\medskip

Comme corollaires des propositions ci-dessus, on obtient les r\'esultats suivants dans le cas o\`u $\mathfrak{q}$ divise $2$ ou $3$. 

\begin{corollary}\label{cor:Phi_2}
On suppose que $\mathfrak{q}$ divise $2$ et que l'une des conditions suivantes est satisfaite~:
\begin{enumerate}
	\item le groupe $\Phi_{\mathfrak{q}}$ est d'ordre  $8$ ou $24$; 
	\item le groupe $\Phi_{\mathfrak{q}}$ est d'ordre  $3$ ou $6$ et le degr\'e r\'esiduel $f_{\mathfrak{q}}$ est impair.
\end{enumerate}
Alors, la repr\'esentation $\rep$ est irr\'eductible pour tout nombre premier $p\geq 5$.
\end{corollary}
Lorsque $\mathfrak{q}$ divise $2$, l'étude faite dans~\cite{Bil09a} permet parfois de calculer l'ordre du groupe~$\Phi_{\mathfrak{q}}$ directement à partir de la valuation de l'invariant modulaire de~$E$ (\cite[th.~1]{Bil09a}). Si $K$ est une extension quadratique de $\Qset$ (ou plus généralement si le degré sur $\Qset_2$ du complété de $K$ en $\mathfrak{q}$ est $\leq 2$), le théorème~\cite[th.~2]{Bil09a} et~\cite{Cal04} fournissent en toute généralité l'ordre du groupe~$\Phi_{\mathfrak{q}}$ en fonction des coefficients d'une équation de Weierstrass de~$E$. 

\medskip

\noindent{\it Remarque. }La condition de parité dans le corollaire précédent est nécessaire. En effet, soient $K$ l'extension de $\Qset$ engendrée par une racine du polynôme
\begin{equation*}
(X^2+5X+1)^3(X^2+13X+49)-\frac{2^4\cdot 13^3}{3^2}X
\end{equation*}
et $E$ la courbe elliptique définie sur~$K$ par l'équation
\begin{equation*}
y^2=x^3-x^2-4x-4.
\end{equation*}
Alors, le degré résiduel de $K$ en l'unique idéal $\mathfrak{p}_2$ de~$\integers$ divisant~$2$, est~$f_{\mathfrak{p}_2}=2$ et la courbe $E$ a un défaut de semi-stabilité d'ordre~$6$ en~$\mathfrak{p}_2$. Pour autant la représentation $\rep[7]: \mathrm{G}_{K}\longrightarrow\mathrm{GL}_2(\Fset_7)$ est \emph{r\'eductible} car $K$ correspond au sous-corps de $\overline{\Qset}$ laissé fixe par le stablilisateur dans $\mathrm{Gal}(\overline{\Qset}/\Qset)$ d'un sous-groupe d'ordre~$7$ de $E(\overline{\Qset})$ (\cite[p.~273]{KrOe92}).

\medskip

Lorsque $\mathfrak{q}$ divise $3$, on a le corollaire suivant.
\begin{corollary}\label{cor:Phi_3}
On suppose que $\mathfrak{q}$ divise $3$ et que l'une des conditions suivantes est satisfaite~:
\begin{enumerate}
	\item le groupe $\Phi_{\mathfrak{q}}$ est d'ordre $12$; 
	\item le groupe $\Phi_{\mathfrak{q}}$ est d'ordre  $4$ et le degr\'e r\'esiduel $f_{\mathfrak{q}}$ est impair.
\end{enumerate}
Alors, la repr\'esentation $\rep$ est irr\'eductible pour tout nombre premier $p\geq 5$.
\end{corollary}

		\section{Rappels}\label{s:rappels}

Dans toute cette section, on fixe un corps de nombres $K$ contenu dans $\Qbar$ et une courbe elliptique $E$ d\'efinie sur $K$. Soit $p$ un nombre premier exceptionnel. Le groupe $\Etors$ possède alors une droite $D$ stable par $\mathrm{G}_{K}$. Notons $\car$ le caractère donnant l'action de $\mathrm{G}_{K}$ sur $D$. On l'appelle caractère d'isogénie associé à $D$. Dans une base convenable de $\Etors$ sur $\Fset_p$, la représentation $\rep$ est représentable matriciellement par
\begin{equation*}
\begin{pmatrix}
 \car & * \\
   0  & \car' \\
\end{pmatrix},
\end{equation*}
où $\car$ et $\car'$ s'interprètent comme des caractères de $\mathrm{G}_{K}$ à valeurs dans $\Fset_p^*$. On a
\begin{equation}\label{eq:determinant_rep}
\det \rep=\car\cdot\car'=\chi_p,
\end{equation}
où $\chi_p$ est le caractère donnant l'action de $\mathrm{G}_K$ sur les racines $p$-ièmes de l'unité (caractère cyclotomique).

La repr\'esentation $\rep$ se factorise \`a travers le groupe de Galois de l'extension $K(\Etors)/K$, o\`u $K(\Etors)$ est le corps engendr\'e sur $K$ par les coordonn\'ees des points de $p$-torsion de $E$. On note encore $\rep$, $\car$, $\car'$ et $\chi_p$ les morphismes pass\'es au quotient. 

Soit $\mathfrak{q}$ est un id\'eal premier de $\integers$. On note $I_{\mathfrak{q}}$ un sous-groupe d'inertie en $\mathfrak{q}$ de $\mathrm{Gal}(K(\Etors)/K)$. Si $E$ a bonne r\'eduction en $\mathfrak{q}$ et $\mathfrak{q}$ ne divise pas $p$, l'extension $K(\Etors)/K$ est non ramifi\'ee en $\mathfrak{q}$ par le crit\`ere de N\'eron-Ogg-Shafarevich. On note $\sigma_{\mathfrak{q}}$ une subsitution de Frobenius en $\mathfrak{q}$ de $\mathrm{Gal}(K(\Etors)/K)$ (bien d\'efinie \`a conjugaison pr\`es). 

Le r\'esultat suivant est bien connu (c.f. \cite[Th.2.4]{Sil86}) et intervient de fa\c con cruciale dans la d\'emonstration des th\'eor\`emes~\ref{th:critere1} et~\ref{th:critere2}.
\begin{proposition}[Hasse~--~Weil]\label{prop:RH}
Soit $\mathfrak{q}$ est un id\'eal premier de $\integers$ en lequel $E$ a bonne r\'eduction. Les racines complexes de $P_{\mathfrak{q}}$ sont de module $N(\mathfrak{q})^{1/2}$. En particulier, on a
\begin{equation*}
\lvert t_{\mathfrak{q}}\rvert\leq 2N(\mathfrak{q})^{1/2}.
\end{equation*}	
Si de plus $\mathfrak{q}$ ne divise pas $p$, le polyn\^ome caract\'eristique de $\rep(\sigma_{\mathfrak{q}})$ est $\overline{P_{\mathfrak{q}}}=P_{\mathfrak{q}}\pmod{p}\in\Fset_p[X]$. En particulier, on a 
\begin{equation*}
\overline{P_{\mathfrak{q}}}(\car(\sigma_{\mathfrak{q}}))=0.
\end{equation*}	
\end{proposition}

		\subsection{L'ensemble $\Exc[E/K]$}

L'objectif de ce \S\ est de démontrer le résultat suivant.
\begin{proposition}\label{prop:CNS_finitude}
Les conditions suivantes sont équivalentes~:
\begin{enumerate}
	\item\label{item:CNS_finitude1} la courbe $E$ n'a pas de multiplication complexe sur $K$ (i.e. $\mathrm{End}_K(E)=\Zset$);
	\item\label{item:CNS_finitude2} l'ensemble $\Exc[E/K]$ est fini.
\end{enumerate}
\end{proposition}
\noindent{\it D\'emonstration. }L'implication $	\ref{item:CNS_finitude1})\Rightarrow \ref{item:CNS_finitude2})$ résulte du théorème de \v{S}afarevi\v{c} sur la finitude des classes de $K$-isomorphisme de courbes elliptiques $K$-isogènes à une courbe donnée (\cite[IX \S6]{Sil86}). Elle est due à Serre et démontrée dans~\cite[IV-9]{Ser68a}.

Réciproquement, si $E$ a des multiplications complexes sur $K$ (i.e. $\mathrm{End}_K(E)$ est de rang~$2$ comme $\Zset$-module), alors 
\begin{equation*}
\mathrm{End}_K(E)\otimes \Qset=\mathrm{End}_{\Qbar}(E)\otimes \Qset 
\end{equation*}
et $K$ contient le corps quadratique imaginaire $L=\mathrm{End}_K(E)\otimes \Qset$. Soit $p$ un nombre premier décomposé dans $L$. On~a
\begin{equation*}
p\integers[L]=\pi\cdot \overline{\pi},
\end{equation*}
où $\integers[L]$ est l'anneau des entiers de $L$. Alors, l'ensemble $E[\pi]$ des points de $E$ annulés par les éléments de $\pi$ est défini sur $K$ et d'ordre~$p$ (\cite[ch.9 \S4]{Lan87}). On en déduit que l'ensemble $\Exc[E/K]$ est infini.

		\subsection{Ramification et caract\`ere d'isog\'enie}\label{ss:ramification_et_caractere}

On suppose que $p$ est un nombre premier exceptionnel pour $E$. Le r\'esultat suivant se d\'eduit de l'\'etude de la restriction de $\rep$ aux sous-groupes d'inertie de $\mathrm{Gal}(K(\Etors)/K)$ telle qu'elle est faite, par exemple, dans~\cite[IV]{Ser68a}, \cite[\S\S1.11-1.12]{Ser72} et~\cite{Kra97a} (voir \'egalement~\cite[\S1]{Dav08} pour une discussion similaire).

\begin{proposition}\label{prop:ramification_caracteres}
Supposons $p\geq 5$ non ramifié dans $K$.
\begin{enumerate}
	\item Le caractère $\car^{12}$ est non ramifié en dehors des idéaux premiers de $\mathcal{O}_K$ divisant $p$.
	\item Soit $\mathfrak{p}$ un idéal de $\integers$ divisant $p$. On suppose que $E$ n'a pas mauvaise réduction additive en~$\mathfrak{p}$ avec potentiellement bonne réduction de hauteur~$2$ (supersingulière). Alors, il existe un entier $\alpha_{\mathfrak{p}}\in\{0,12\}$ tel que 
\begin{equation*}
	\car^{12}\mid_{I_{\mathfrak{p}}}=\chi_p^{\alpha_{\mathfrak{p}}}\mid_{I_{\mathfrak{p}}}.
\end{equation*}
\end{enumerate}
\end{proposition}

\noindent{\it Remarques. }\label{rem:convention}
\begin{enumerate}
	\item Dans une base convenable, la représentation sur les points de $p$-torsion de la courbe $E/D$ est représentable matriciellement par
\begin{equation*}
\begin{pmatrix}
 \car' & * \\
   0  & \car \\
\end{pmatrix}.
\end{equation*}
Autrement dit, d'après l'égalité~(\ref{eq:determinant_rep}), on peut toujours, si on le souhaite, remplacer la famille $\{\alpha_{\mathfrak{p}}\}_{\mathfrak{p}\mid p}$ par la famille $\{12-\alpha_{\mathfrak{p}}\}_{\mathfrak{p}\mid p}$.
	\item On peut montrer en utilisant la description locale de $\rep$ donnée dans la proposition \cite[prop.$2$]{Kra97a} que si $\mathfrak{p}$ divise $p$ et $E$ a mauvaise r\'eduction additive en $\mathfrak{p}$ avec potentiellement bonne r\'eduction supersinguli\`ere, alors il existe un entier $\alpha_{\mathfrak{p}}\in\{4,6,8\}$ tel que 
\begin{equation*}
	\car^{12}\mid_{I_{\mathfrak{p}}}=\chi_p^{\alpha_{\mathfrak{p}}}\mid_{I_{\mathfrak{p}}}.
\end{equation*}
	\item Dans sa th\`ese (\cite{Dav08}), A. David d\'emontre que si $K$ ne contient pas le corps de classes de Hilbert d'un corps quadratique imaginaire, il existe alors une constante effective $C(K)$, ne d\'ependant que de $K$ (et donc pas de $E$) telle que si $p>C(K)$, on a $\alpha_{\mathfrak{p}}=6$ pour \emph{tout} id\'eal premier $\mathfrak{p}$ de $\integers$ divisant~$p$ (voir \'egalement \cite{Mom95}). Nous n'utiliserons pas ces r\'esultats.
\end{enumerate}

	\subsection{Th\'eorie du corps de classes et caract\`ere d'isog\'enie}\label{ss:lemme_Momose}

On reprend les hypoth\`eses et notations pr\'ec\'edentes. En particulier, $p$ est un nombre premier $\geq 5$ non ramifié dans $K$ et on suppose que pour tout idéal premier $\mathfrak{p}$ de~$\integers$ divisant $p$, $E$ n'a pas mauvaise réduction additive en~$\mathfrak{p}$ avec potentiellement bonne réduction de hauteur~$2$. \'Etant donné un idéal premier $\mathfrak{p}$ de $\integers$ au-dessus de $p$, on désigne par
\begin{equation*}
N_{\mathfrak{p}}:(\integers/\mathfrak{p})^*\longrightarrow \Fset_p^*
\end{equation*}
le morphisme norme. L'objectif de ce paragraphe~\ref{ss:lemme_Momose} est de démontrer la proposition ci-dessous, cruciale dans la d\'emonstration des th.~\ref{th:critere1} et~\ref{th:critere2}. Elle figure également sous une forme légèrement différente dans la thèse de David (\cite[prop.~2.2.1]{Dav08}) ainsi que dans l'article~\cite[lem.~1]{Mom95} de Momose (sous l'hypoth\`ese que $K/\Qset$ est galoisienne).
\begin{proposition}\label{prop:CFT}
Soit $a\in\integers$ premier \`a~$p$ et $a\integers=\prod_{\mathfrak{q}}\mathfrak{q}^{v_{\mathfrak{q}}(a)}$ la d\'ecomposition de $a\integers$ en produit d'id\'eaux premiers de~$\integers$. On suppose que pour tout idéal premier $\mathfrak{q}$ de $\integers$ divisant $a$, $E$ a bonne réduction en $\mathfrak{q}$. Alors, on~a~:
\begin{equation*}
\prod_{\mathfrak{q}\mid a}\car(\sigma_{\mathfrak{q}})^{12v_{\mathfrak{q}}(a)}=\prod_{\mathfrak{p}\mid p}N_{\mathfrak{p}}(a+\mathfrak{p})^{\alpha_{\mathfrak{p}}},
\end{equation*}
o\`u $\alpha_{\mathfrak{p}}\in\{0,12\}$ est d\'efini \`a la proposition~\ref{prop:ramification_caracteres}.
\end{proposition}

	\subsubsection{Un lemme de la théorie du corps de classes.}

Soient $L$ l'extension de $K$ trivialisant le caract\`ere $\car^{12}$ et $\mu_p$ le groupe de racines $p$-ièmes de l'unité dans~$\Qbar$. D'après l'accouplement de Weil, on a $\mu_p\subset K(\Etors)$. Donc $L(\mu_p)$ est une sous-extension ab\'elienne de $K(\Etors)/K$. On note $\ide$ le groupe des id\`eles de $K$ et
\begin{equation*}
\rec{}~:\ide\longrightarrow \mathrm{Gal}(L(\mu_p)/K),
\end{equation*}
le morphisme de r\'eciprocit\'e global donn\'e par la th\'eorie du corps de classes. Il est surjectif et son noyau contient les id\`eles principales.

Soit $v$ une place de $K$. On note $K_v$ le compl\'et\'e de $K$ en $v$ et on identifie $K$ à un sous-corps de $K_v$. On désigne par
\begin{equation*}
\rec{v}~:K_v^*\hookrightarrow\ide\longrightarrow \mathrm{Gal}(L(\mu_p)/K)
\end{equation*}
la compos\'ee de l'injection de $K_v^*$ dans $\ide$ par le morphisme de r\'eciprocit\'e global.

Si $\mathfrak{q}$ est un id\'eal premier de $\integers$ de bonne r\'eduction ne divisant pas $p$, on rappelle que l'extension $K(\Etors)/K$ est non ramifi\'ee en $\mathfrak{q}$. La restriction \`a $\mathrm{Gal}(L(\mu_p)/K)$ d'une substitution de Frobenius en~$\mathfrak{q}$ du groupe $\mathrm{Gal}(K(\Etors)/K)$ (bien d\'efinie \`a conjugaison pr\`es) est unique. On la note encore~$\sigma_{\mathfrak{q}}$. De m\^eme, on note encore $\chi_p$ (resp. $\car$) la restriction du caract\`ere cyclotomique (resp. d'isog\'enie) \`a $\mathrm{Gal}(L(\mu_p)/K)$.

Le lemme suivant regroupe plusieurs résultats classiques de la théorie du corps de classes qui seront utiles à la démonstration de la proposition~\ref{prop:CFT}. La démonstration du troisième point est tirée de \cite[App.~1 prop.~1]{Kra07}.

\begin{lemma}\label{lem:CFT} Soit $v$ une place de $K$.
\begin{enumerate}
	\item Si $v$ est une place infinie de $K$, on a $\car^{12}(r_v(a))=1$.
	\item Si $v=\mathfrak{q}$ est une place finie de $K$ ne divisant pas $p$, on a $\rec{\mathfrak{q}}(\mathcal{U}_{\mathfrak{q}})=\{1\}$, o\`u $\mathcal{U}_{\mathfrak{q}}$ est le groupe des unit\'es de l'anneau d'entiers du corps $K_\mathfrak{q}$. Si de plus, $\mathfrak{q}$ divise~$a$, alors $\rec{\mathfrak{q}}(\pi_\mathfrak{q})=\sigma_{\mathfrak{q}}$, où $\pi_\mathfrak{q}$ est une uniformisante de $K_\mathfrak{q}$.
	\item Si $v=\mathfrak{p}$ est une place finie de $K$ divisant $p$, alors $\rec{\mathfrak{p}}(a)$ appartient au sous-groupe d'inertie en~$\mathfrak{p}$ de $L(\mu_p)/K$ et on a 
\begin{equation*}
\chi_p \left(\rec{\mathfrak{p}}(a)\right)=N_{\mathfrak{p}}(a+\mathfrak{p})^{-1}.
\end{equation*}
\end{enumerate}
\end{lemma}
\noindent{\it D\'emonstration. }Soit $v$ une place de $K$. On distingue trois cas.
\begin{enumerate}
	\item Supposons que $v$ soit une place infinie de $K$. Soit $L'$ l'extension de $K$ trivialisant le caract\`ere $\car$, 
\begin{equation*}
\rec{}'~:\ide\longrightarrow \mathrm{Gal}(L'(\mu_p)/K),
\end{equation*}
le morphisme de r\'eciprocit\'e global donn\'e par la th\'eorie du corps de classes et 
\begin{equation*}
\rec{v}'~:K_v^*\hookrightarrow\ide\stackrel{r'}{\longrightarrow} \mathrm{Gal}(L'(\mu_p)/K).
\end{equation*}
L'image de l'application $\rec{v}'$ est d'ordre $\leq 2$. Par ailleurs, l'image par $\car^{12}$ d'un \'el\'ement de $\mathrm{Gal}(L'(\mu_p)/K)$ ne d\'epend que de sa restriction \`a $\mathrm{Gal}(L(\mu_p)/K)$. D'o\`u~:
\begin{equation*}
\car^{12}(\rec{v}(a))=\car^{12}(\rec{v}'(a)),
\end{equation*}
puis
\begin{equation*}
\car(\rec{v}'(a))^{12}=\car(\rec{v}'(a)^{12})=1.
\end{equation*}
D'où le résultat.

	\item Supposons que $v=\mathfrak{q}$ soit une place finie de $K$ ne divisant pas $p$. Alors, d'après~\cite{Neu86}, l'image par $\rec{\mathfrak{q}}$ de $\mathcal{U}_{\mathfrak{q}}$ est un sous-groupe d'inertie en $\mathfrak{q}$ de l'extension $L(\mu_p)/K$. Or celle-ci est non ramifiée en $\mathfrak{q}$ d'après le crit\`ere de N\'eron-Ogg-\v{S}afarevi\v{c}. D'où l'égalité 
\begin{equation*}
\rec{\mathfrak{q}}(\mathcal{U}_{\mathfrak{q}})=\{1\}.
\end{equation*}
Si de plus $\mathfrak{q}$ divise~$a$ alors $E$ a bonne réduction en $\mathfrak{q}$ par hypoth\`ese et d'après~\cite{Neu86}, l'image par $\rec{\mathfrak{q}}$ de $\pi_{\mathfrak{q}}$ est la substitution de Frobenius en $\mathfrak{q}$ de l'extension $L(\mu_p)/K$. Autrement dit, $\rec{\mathfrak{q}}(\pi_\mathfrak{q})=\sigma_{\mathfrak{q}}$.

	\item Supposons que $v=\mathfrak{p}$ soit une place finie de $K$ divisant $p$. On note $\overline{\Qset_p}$ une clôture algébrique de $\Qset_p$. Comme $p$ est non ramifié dans $K$, on identifie $K_{\mathfrak{p}}$ à l'extension non ramifiée de $\Qset_p$ contenue dans $\overline{\Qset_p}$ dont le degré sur $\Qset_p$ est le degré résiduel de $\mathfrak{p}$ sur~$p$. On note $K^{ab}$ la clôture abélienne de $K$ dans $\Qbar$, $K_{\mathfrak{p}}^{ab}$ la clôture abélienne de $K_{\mathfrak{p}}$ dans $\overline{\Qset_p}$,
\begin{equation*}
\Theta_{\mathfrak{p}}~: K_{\mathfrak{p}}^*\longrightarrow \mathrm{Gal}(K_{\mathfrak{p}}^{ab}/K_{\mathfrak{p}})
\end{equation*}
le morphisme de réciprocité local en~${\mathfrak{p}}$ et
\begin{equation*}
\mathrm{Res}_{\mathfrak{p}}~:\mathrm{Gal}(K_{\mathfrak{p}}^{ab}/K_{\mathfrak{p}})\longrightarrow \mathrm{Gal}(L(\mu_p)/K)
\end{equation*}
le morphisme de restriction. D'après la compatibilité entre la théorie du corps de classes locale et globale, on a, pour tout $x\in K_{\mathfrak{p}}^*$,
\begin{equation}\label{eq:compa_local_global}
\mathrm{Res}_{\mathfrak{p}}(\Theta_{\mathfrak{p}}(x))=r_{\mathfrak{p}}(x).
\end{equation}

Or, d'après le corollaire de~\cite[App.~1 prop.~1]{Kra07}, on a 
\begin{equation*}
\Theta_{\mathfrak{p}}(a)(\zeta)=\zeta^{n^{-1}},
\end{equation*}
où $\zeta$ est une racine primitive $p$-i\`eme de l'unit\'e dans~$\overline{\Qset_p}$ et $n$ est un entier tel que 
\begin{equation*}
N_{\mathfrak{p}}(a+\mathfrak{p})=n\pmod{p\Zset}.
\end{equation*}
D'où le résultat voulu, d'après l'égalité~(\ref{eq:compa_local_global}). 
\end{enumerate}
Cela termine la démonstration du lemme~\ref{lem:CFT}.

		\subsubsection{Démonstration de la proposition~\ref{prop:CFT}.}

L'entier $a$ est non nul car premier \`a~$p$. L'image par le morphisme de r\'eciprocit\'e global de l'id\`ele principale $(a)_v$ est triviale~: 
\begin{equation}\label{eq:CFT_produit}
\prod_{v}\rec{v}(a)=1.
\end{equation}
Si $v$ est une place infinie de $K$, alors d'apr\`es le lemme~\ref{lem:CFT}, on a
\begin{equation}\label{eq:CFT_infini}
\car^{12}\left(\rec{v}(a)\right)=1.
\end{equation}
Si $v=\mathfrak{q}$ est une place finie de $K$ ne divisant ni $p$, ni $a$, alors $a\in \mathcal{U}_{K_\mathfrak{q}}$. D'apr\`es loc.~cit. on a donc $\rec{\mathfrak{q}}(a)=1$. 

Si $v=\mathfrak{q}$ est une place finie de $K$ divisant~$a$. Alors, $a=u\cdot \pi_{\mathfrak{q}}^{v_{\mathfrak{q}}(a)}$,
o\`u $\pi_{\mathfrak{q}}$ est une uniformisante de $K_\mathfrak{q}$ et $u\in\mathcal{U}_{K_\mathfrak{q}}$. D'apr\`es loc.~cit. on a donc $\rec{\mathfrak{q}}(a)=\sigma_{\mathfrak{q}}^{v_{\mathfrak{q}}(a)}$, puis
\begin{equation}\label{eq:CFT_hors_p}
\car^{12}\left(\rec{\mathfrak{q}}(a)\right)=(\car^{12}(\sigma_{\mathfrak{q}}))^{v_{\mathfrak{q}}(a)}=\car(\sigma_{\mathfrak{q}})^{12v_{\mathfrak{q}}(a)}.
\end{equation}
Si $v=\mathfrak{p}$ est une place finie de $K$ divisant $p$, alors d'apr\`es loc.~cit., $\rec{\mathfrak{p}}(a)$ appartient au sous-groupe d'inertie en~$\mathfrak{p}$ de $L(\mu_p)/K$ et on a
\begin{equation*}
\chi_p \left(\rec{\mathfrak{p}}(a)\right)=N_{\mathfrak{p}}(a+\mathfrak{p})^{-1}.
\end{equation*}
Or, d'apr\`es la proposition~\ref{prop:ramification_caracteres}, on a
\begin{equation*}
	\car^{12}\mid_{I_{\mathfrak{p}}}=\chi_p^{\alpha_{\mathfrak{p}}}\mid_{I_{\mathfrak{p}}}.
\end{equation*}
On en d\'eduit que l'on a
\begin{equation}\label{eq:CFT_en_p}
	\car^{12}\left(\rec{\mathfrak{p}}(a)\right)=N_{\mathfrak{p}}(a+\mathfrak{p})^{-\alpha_{\mathfrak{p}}}.
\end{equation}
D'apr\`es les \'egalit\'es~(\ref{eq:CFT_produit})--(\ref{eq:CFT_en_p}) ci-dessus, on a
\begin{align*}
1	&=\prod_{v}\car^{12}\left(\rec{v}(a)\right) \\
	&=\prod_{\mathfrak{q}\mid a}\car(\sigma_{\mathfrak{q}})^{12v_{\mathfrak{q}}(a)}\cdot\prod_{\mathfrak{p}\mid p}N_{\mathfrak{p}}(a+\mathfrak{p})^{-\alpha_{\mathfrak{p}}}.
\end{align*}
Cela d\'emontre la proposition~\ref{prop:CFT}.

	\section{D\'emonstrations des th\'eor\`emes~\ref{th:critere1} et~\ref{th:critere2}}

Soient $K$ un corps de nombres contenu dans $\Qbar$, $E$ une courbe elliptique d\'efinie sur $K$ et $p$ un nombre premier exceptionnel pour~$(E,K)$. On suppose $p\geq 5$, non ramifié dans~$K$ et pour tout idéal premier $\mathfrak{p}$ de $\integers$ divisant $p$, $E$ n'a pas mauvaise réduction additive en~$\mathfrak{p}$ avec potentiellement bonne réduction de hauteur~$2$.

Avant de d\'emontrer les th\'eor\`emes~\ref{th:critere1} et~\ref{th:critere2}, on commence par d\'efinir pour tout anneau intègre~$A$, une loi de mono\"ide commutatif~$\ast$ sur un sous-ensemble de $A[X]$ et par en \'etudier les propri\'et\'es utiles.

	\subsection{Loi de monoïde}\label{ss:lemmes_preliminaires}

Soit $A$ un anneau intègre de corps des fractions $L$ et $\overline{L}$ une clôture algébrique de~$L$. On note $M_A$ le sous-ensemble de $A[X]$ constitué des polynômes unitaires ne s'annulant pas en~$0$.
\begin{lemma}\label{lem:loi_monoide}
L'application
\begin{equation*}
\begin{array}{ccl}
M_A\times M_A & \longrightarrow & A[X] \\
(P,Q) & \longmapsto & (P\ast Q)(X)=\Res_Z\left(P(Z),Q\left(X/Z\right)Z^{\deg Q}\right)
\end{array}
\end{equation*}
a une image contenue dans $M_A$. Elle définit une loi de monoïde commutatif sur~$M_A$ d'élément neutre~$\Psi_1(X)=X-1$. De plus, si $P$, $Q\in M_A$ s'écrivent 
\begin{equation*}
P(X)=\prod_{i=1}^n(X-\alpha_i)\quad\textrm{et}\quad Q(X)=\prod_{j=1}^m(X-\beta_j)
\end{equation*}
dans $\overline{L}[X]$, on a 
\begin{equation*}
(P\ast Q)(X)=\prod_{\substack{1\leq i\leq n\\ 1\leq j\leq m}}(X-\alpha_i\beta_j).
\end{equation*}
En particulier,
\begin{equation*}
(P\ast Q)(0)=(-1)^{\deg P\cdot \deg Q}P(0)^{\deg Q}Q(0)^{\deg P}.
\end{equation*}

\end{lemma}
\noindent{\it D\'emonstration. }Il s'agit de vérifier que pour tout $P$, $Q$ et $R\in M_A$, on a 
\begin{enumerate}
	\item $P\ast Q\in M_A$;
	\item $P\ast \Psi_1=\Psi_1\ast P=P$;
	\item $(P\ast Q)\ast R=P\ast (Q\ast R)$;
	\item $P\ast Q=Q\ast P$.
\end{enumerate}
On suppose que le polynôme $Q$ s'écrit
\begin{equation*}
Q(X)=X^m+b_{m-1}X^{m-1}+\cdots+b_1X+b_0,\quad\textrm{avec }b_0\not=0.
\end{equation*}
Alors, 
\begin{equation*}
Q\left(\frac{X}{Z}\right)Z^m=b_0Z^m+b_{1}XZ^{m-1}+\cdots+b_{m-1}X^{m-1}Z+X^m\in A[X][Z]
\end{equation*}
et $\deg_Z(Q(X/Z)Z^m)=m=\deg Q$ (car $b_0\not=0$). Par définition du résultant de deux polynômes (\cite[A IV.71 \S6 Déf.~1]{Bou81}), on a donc $P\ast Q\in A[X]$. Par ailleurs, sur $\overline{L}$, on a
\begin{equation*}
Q\left(\frac{X}{Z}\right)Z^m=Q(0)\prod_{j=1}^m\left(Z-\frac{1}{\beta_j}X\right)
\end{equation*}
et d'après~\cite[A IV.75 \S6 Cor.~1]{Bou81},
\begin{equation*}
(P\ast Q)(X)=Q(0)^n\prod_{i,j}\left(\alpha_i-\frac{1}{\beta_j}X\right).
\end{equation*}
Or, $Q(0)=\prod_{j=1}^m(-\beta_j)$, donc
\begin{equation*}
(P\ast Q)(X)=\prod_{i=1}^n\prod_{j=1}^m\left(X-\alpha_i\beta_j\right).
\end{equation*}
C'est la formule de l'énoncé. On en déduit que l'on a~:
\begin{itemize}
	\item $(P\ast Q)(0)=(-1)^{\deg P\cdot \deg Q}P(0)^{\deg Q}Q(0)^{\deg P}\not= 0$, donc $P\ast Q\in M_A$;
	\item $P\ast \Psi_1=\Psi_1\ast P=P$;
	\item $P\ast Q=Q\ast P$;
\end{itemize}
De plus, les polynômes $(P\ast Q)\ast R$ et $P\ast (Q\ast R)$ ont les mêmes racines dans $\overline{L}$ comptées avec multiplicités. Comme ils sont unitaires, ils sont égaux. D'où le lemme.

\bigskip

\begin{lemma}\label{lem:morphisme_monoide}
Soient $r\geq 1$ et $P\in M_A$. Il existe un unique polynôme $P^{(r)}\in M_A$ tel que
\begin{equation}\label{eq:image_morphisme_monoide}
P^{(r)}(X^r)=(P\ast \Psi_r)(X)
\end{equation}
où $\Psi_r(X)=X^r-1$. L'application $P\mapsto P^{(r)}$ est un morphisme de monoïdes pour la loi~$\ast$. De plus, si $P\in M_A$ se factorise sur $\overline{L}$ de la façon suivante
\begin{equation}\label{eq:expression_de_P(r)}
P(X)=\prod_{i=1}^n(X-\alpha_i),\quad\textrm{on a}\quad P^{(r)}(X)=\prod_{i=1}^n\left(X-\alpha_i^r\right).
\end{equation}
\end{lemma}
\noindent{\it D\'emonstration. }Soit $P\in M_A$. L'unicité d'un polynôme~$P^{(r)}$ vérifiant la relation~(\ref{eq:image_morphisme_monoide}) est immédiate. Posons
\begin{equation*}
P(X)=\prod_{i=1}^n(X-\alpha_i)\quad\textrm{avec } \alpha_i\in\overline{L}
\end{equation*}
et $\zeta_r$ une racine $r$-ième de l'unité dans~$\overline{L}$. Vérifions alors qu'il existe bien un polynôme de $M_A$ satisfaisant à l'égalité~(\ref{eq:image_morphisme_monoide})~:
\begin{align*}
(P\ast \Psi_r)(X) 	& = \Res_Z\left(P(Z),X^r-Z^r\right) \\
		& =(-1)^n\prod_{i=1}^n\prod_{k=0}^{r-1}\left(\alpha_i-\zeta_r^kX\right)\\
		& \qquad\textrm{d'après  }\cite[\textrm{A IV.75 \S6 Cor.~1}]{Bou81}\\
 		& =(-1)^n \prod_{i=1}^n(-1)^r\zeta_r^{\frac{r(r-1)}{2}}\prod_{k=0}^{r-1}\left(X-\zeta_r^{-k}\alpha_i\right)
\end{align*}
Or, on a $\zeta_r^{\frac{r(r-1)}{2}}=(-1)^{r+1}$. D'où
\begin{equation*}
(P\ast \Psi_r)(X)=\prod_{i=1}^n\prod_{k=0}^{r-1}\left(X-\zeta_r^{-k}\alpha_i\right)=\prod_{i=1}^n\left(X^r-\alpha_i^r\right).
\end{equation*}
Autrement dit, 
\begin{equation*}
(P\ast \Psi_r)(X)=P^{(r)}(X^r)
\end{equation*}
où l'on a posé
\begin{equation*}
P^{(r)}(X)=\prod_{i=1}^n\left(X-\alpha_i^r\right).
\end{equation*}
Cela démontre qu'il existe bien un polynôme $P^{(r)}$ de $A[X]$ satisfaisant à l'égalité~(\ref{eq:image_morphisme_monoide}) et qu'il est donné par la formule~(\ref{eq:expression_de_P(r)}). Par ailleurs, d'après le lemme~\ref{lem:loi_monoide}, on a $P^{(r)}(0)=(-1)^{(r+1)\deg P}P(0)^r\not=0$ et comme $P^{(r)}$ est unitaire, on a $P^{(r)}\in M_A$. On en déduit que l'application $P\mapsto P^{(r)}$ est bien définie. 

Vérifions enfin qu'il s'agit bien d'un morphisme de monoïdes. On a $\Psi_1^{(r)}=\Psi_1$. Soient $P$ et $Q$ dans $M_A$. D'après le lemme~\ref{lem:loi_monoide} et la formule~(\ref{eq:expression_de_P(r)}), les polynômes $(P\ast Q)^{(r)}$ et $P^{(r)}\ast Q^{(r)}$ ont les mêmes racines dans $\overline{L}$ comptées avec multiplicités. Ils sont donc égaux. D'où le lemme~\ref{lem:morphisme_monoide}.

\bigskip

\begin{lemma}\label{lem:morphisme_anneaux_morphisme_monoides}
Soient $A$ et $B$ deux anneaux intègres et ${\varphi~:A\rightarrow B}$ un morphisme d'anneaux. L'ensemble 
\begin{equation*}
M_A^{\varphi}=\{P\in M_A\mid \varphi(P(0))\not=0\}
\end{equation*}
est stable pour la loi $\ast$. L'application $\varphi$ induit un morphisme de monoïdes (encore noté~$\varphi$)
\begin{equation*}
\varphi~:M_A^{\varphi}\longrightarrow M_B.
\end{equation*}
Soient $P\in M_A^{\varphi}$ et $r\geq 1$. Alors, $P^{(r)}\in M_A^{\varphi}$ et on a $\left(\varphi(P)\right)^{(r)}=\varphi\big(P^{(r)}\big)$.
\end{lemma}
\noindent{\it D\'emonstration. }D'après le lemme~\ref{lem:loi_monoide}, si $P$, $Q\in M_A^{\varphi}\subset M_A$, on a $P\ast Q\in M_A$ et 
\begin{equation*}
(P\ast Q)(0)=(-1)^{\deg P\cdot \deg Q}P(0)^{\deg Q}Q(0)^{\deg P}.
\end{equation*}
D'où 
\begin{equation*}
\varphi\left((P\ast Q)(0)\right)=(-1)^{\deg P\cdot \deg Q}\varphi(P(0))^{\deg Q}\varphi(Q(0))^{\deg P}\not=0
\end{equation*}
car $\varphi(P(0))\not=0$, $\varphi(Q(0))\not=0$ et $B$ est intègre. Donc $M_A^{\varphi}$ est bien un sous-ensemble de $M_A$ stable pour la loi~$\ast$. On a $\varphi(\Psi_1)=\Psi_1$ et la relation
\begin{equation*}
\varphi(P\ast Q)=\varphi(P)\ast \varphi(Q),\quad\textrm{pour }P,\textrm{ }Q\in M_A^{\varphi},
\end{equation*}
résulte de la définition du résultant de deux polynômes (\cite[A IV.72 \S6]{Bou81}) en termes de déterminant de Sylvester.

Soient $P\in M_A^{\varphi}$ et $r\geq 1$. Alors, d'après le lemme~\ref{lem:morphisme_monoide}, on a $P^{(r)}\in M_A$ et
\begin{equation*}
\varphi\big(P^{(r)}(0)\big)=(-1)^{(r+1)\deg P}\varphi(P(0))^r\not=0
\end{equation*}
car $\varphi(P(0))\not=0$ et $B$ est intègre. D'où $P^{(r)}\in M_A^{\varphi}$. De plus, d'après la formule~(\ref{eq:image_morphisme_monoide}) et la définition du résultant de deux polynômes (\cite[A IV.72 \S6]{Bou81}), on a
\begin{equation*}
\varphi\big(P^{(r)}(X^r)\big)=\varphi\big((P\ast \Psi_r)(X)\big)=\big(\varphi(P)\ast \Psi_r)(X)=\varphi(P)^{(r)}(X^r).
\end{equation*}
D'où l'égalité $\left(\varphi(P)\right)^{(r)}=\varphi\big(P^{(r)}\big)$ et le lemme~\ref{lem:morphisme_anneaux_morphisme_monoides}.

\medskip

\noindent{\it Remarque. }D'après le lemme ci-dessus, l'application de réduction $\varphi~:\Zset\rightarrow \Zset/p\Zset$ induit un mor\-phisme de monoïdes
\begin{equation*}
\begin{array}{ccl}
M_{\Zset}^{\varphi} & \longrightarrow & M_{\Fset_p} \\
P & \longmapsto & \overline{P}.
\end{array}
\end{equation*}
En particulier, $\overline{P\ast Q}=\overline{P}\ast \overline{Q}$ pour tout $P$, $Q\in M_{\Zset}^{\varphi}$.

	\subsection{D\'emonstration du th\'eor\`eme~\ref{th:critere1}}

Soit $\ell$ un nombre premier tel que $E$ ait bonne r\'eduction en tout id\'eal premier de $\integers$ divisant~$\ell$. Il s'agit de montrer que $p$ divise~$B_{\ell}$.

	\subsubsection{Le polyn\^ome $P_{\ell}^*$}

Soit $\prod_{\ideal\mid\ell}\ideal^{v_{\ideal}(\ell)}$ la décomposition de l'idéal $\ell\integers$ en produit d'idéaux premiers de~$\integers$. On note $g_{\ell}$ le cardinal de l'ensemble $\{\ideal\mid\ell\}$. On rappelle que $E$ a par hypothèse bonne réduction en tout idéal premier $\ideal$ de $\integers$ divisant~$\ell$. Le polynôme $P_{\ell}^*$ donné par la formule~(\ref{eq:polynôme_P_ell*}) est alors bien défini et à coefficients entiers. 
\begin{lemma}\label{lem:polynome_P_ell}
Le polynôme $P_{\ell}^*$ appartient à $M_{\Zset}$ et vérifie~:
\begin{equation}\label{eq:P_ell*_en_0}
P_{\ell}^*(0)=\ell^{12\cdot d\cdot 2^{g_{\ell}-1}}.
\end{equation}
Ses racines complexes sont de module~$\ell^{6d}$. Si de plus $\ell\not=p$, alors $P_{\ell}^*\in M_{\Zset}^{\varphi}$ et on a
\begin{equation*}
\overline{P_{\ell}^*}(\Omega)=0,\quad\textrm{où}\quad \Omega=\prod_{\ideal\mid\ell}\car(\sigma_{\ideal})^{12v_{\ideal}(\ell)}\in\Fset_p.
\end{equation*}
\end{lemma}
\noindent{\it D\'emonstration. }Pour tout $\ideal\mid\ell$, le polynôme $P_{\ideal}$ est unitaire, à coefficients entiers et on a (prop.~\ref{prop:RH})~:
\begin{equation*}
P_{\ideal}(0)=N(\ideal)=\ell^{f_{\ideal}}.
\end{equation*} 
En particulier, $P_{\ideal}\in M_{\Zset}$. D'après les lemmes~\ref{lem:loi_monoide} et~\ref{lem:morphisme_monoide}, le polynôme $P_{\ell}^*$ est bien défini (la loi $\ast$ est associative) et indépendant de l'ordre des idéaux premiers dans la décomposition de~$\ell$ dans~$K$ (la loi $\ast$ est commutative). De plus, $P_{\ell}^*$ appartient à $M_{\Zset}\subset\Zset[X]$. 

Soient $P_1,\ldots,P_n\in M_A$ de degrés respectifs $d_1,\ldots,d_n$. On montre par récurrence sur $n$, à partir de la formule pour $n=2$ du lemme~\ref{lem:loi_monoide} que l'on a
\begin{equation*}
(P_1\ast \cdots \ast P_n)(0)=(-1)^{(n+1)d_1\cdots d_n}\prod_{i=1}^nP_i(0)^{\prod_{j\not=i}d_j}.
\end{equation*}
De plus, d'après le lemme~\ref{lem:morphisme_monoide}, pour tout $P\in M_{\Zset}$ et tout entier $r\geq 1$, on a 
\begin{equation*}
P^{(r)}(0)=(-1)^{(r+1)\deg P}P(0)^r.
\end{equation*}
Comme pour tout idéal $\ideal\mid\ell$, on a $\deg P_{\ideal}=2$, on en déduit
\begin{equation*}
P_{\ell}^*(0)=\prod_{\ideal\mid\ell}P_{\ideal}(0)^{12v_{\ideal}(\ell)\cdot 2^{g_{\ell}-1}}=\prod_{\ideal\mid\ell}\left(\ell^{f_{\ideal}}\right)^{12v_{\ideal}(\ell)\cdot 2^{g_{\ell}-1}}=\ell^{12\cdot 2^{g_{\ell}-1}\sum_{\ideal\mid\ell}f_{\ideal}v_{\ideal}(\ell)}.
\end{equation*}
D'où la formule car $\sum_{\ideal\mid\ell}f_{\ideal}v_{\ideal}(\ell)=d$. 

Par ailleurs, d'après la proposition~\ref{prop:RH}, les racines complexes de $P_{\ideal}$ sont de module~$N(\ideal)^{1/2}=\ell^{f_{\ideal}/2}$. Donc, d'après le lemme~\ref{lem:morphisme_monoide}, celles de $P_{\ideal}^{(12v_{\ideal}(\ell))}$ sont de module~$\ell^{6f_{\ideal}v_{\ideal}(\ell)}$. D'après le lemme~\ref{lem:loi_monoide}, celles de $P_{\ell}^*$ sont de module
\begin{equation*}
\prod_{\ideal\mid\ell}\ell^{6f_{\ideal}v_{\ideal}(\ell)}=\ell^{6\sum_{\ideal\mid\ell}f_{\ideal}v_{\ideal}(\ell)}=\ell^{6d}.
\end{equation*}
Supposons à présent $\ell\not=p$. Alors, d'après la formule~(\ref{eq:P_ell*_en_0}), on a $P_{\ell}^*\in M_{\Zset}^{\varphi}$. D'après la proposition~\ref{prop:RH}, on a
\begin{equation*}
\overline{P_{\ideal}}(\car(\sigma_{\ideal}))=0.
\end{equation*}
Donc d'après le lemme~\ref{lem:morphisme_monoide}, on~a~:
\begin{equation}\label{eq:racines_mod_p_de_P_qd}
\overline{P_{\ideal}}^{(12v_{\ideal}(\ell))}\left(\car(\sigma_{\ideal})^{12v_{\ideal}(\ell)}\right)=0\pmod{p}.
\end{equation}
Puis,
\begin{align*}
\overline{P_{\ell}^*}(\Omega)	& = \overline{\myprod_{\ideal\mid\ell}P_{\ideal}^{(12v_{\ideal}(\ell))}}\Big(\prod_{\ideal\mid\ell}\car(\sigma_{\ideal})^{12v_{\ideal}(\ell)}\Big) \\
				& = \left(\myprod_{\ideal\mid\ell}\overline{P_{\ideal}}^{(12v_{\ideal}(\ell))}\right)\Big(\prod_{\ideal\mid\ell}\car(\sigma_{\ideal})^{12v_{\ideal}(\ell)}\Big)\quad\textrm{(lemme~\ref{lem:morphisme_anneaux_morphisme_monoides})}\\
				& = 0\pmod{p}\quad\textrm{(d'après le lemme~\ref{lem:loi_monoide} et la relation~(\ref{eq:racines_mod_p_de_P_qd}))}.
\end{align*}
D'où le lemme~\ref{lem:polynome_P_ell}.

	\subsubsection{Fin de la démonstration du théorème~\ref{th:critere1}}\label{sss:fin_demo_th1}

Supposons $p=\ell$. Alors, pour $d\geq 2$, par définition de $B_{p}^{(d)}$, il existe un entier $k>0$ tel que $P_p^*(p^{12k})$ divise $B_{p}^{(d)}$. D'où $p$ divise $B_p^{(d)}$ car d'après le lemme~\ref{lem:polynome_P_ell}~:
\begin{equation*}
P_p^*(p^{12k})\equiv P_p^*(0)\equiv 0\pmod{p}.
\end{equation*}

Supposons $p\not=\ell$. D'après la proposition~\ref{prop:CFT} appliqu\'ee \`a $a=\ell$, on a~:
\begin{equation}\label{eq:demo_egalite_CFT}
\Omega=\prod_{\ideal\mid\ell}\car(\sigma_{\ideal})^{12v_{\ideal}(\ell)}=\prod_{\mathfrak{p}\mid p}N_{\mathfrak{p}}(\ell+\mathfrak{p})^{\alpha_{\mathfrak{p}}}.
\end{equation}
Or, par définition on a 
\begin{equation*}
N_{\mathfrak{p}}(\ell+\mathfrak{p})=\ell^{1+p+\cdots+p^{f_{\mathfrak{p}}-1}}\pmod{p}=\ell^{f_{\mathfrak{p}}}\pmod{p}
\end{equation*}
où $N(\mathfrak{p})=\lvert\integers/\mathfrak{p}\rvert=\ell^{f_{\mathfrak{p}}}$. D'où
\begin{equation}\label{eq:RHS}
\prod_{\mathfrak{p}\mid p}N_{\mathfrak{p}}(\ell+\mathfrak{p})^{\alpha_{\mathfrak{p}}}=\ell^{\sum_{\mathfrak{p}\mid p}f_{\mathfrak{p}}\alpha_{\mathfrak{p}}}.
\end{equation}
Or, $\alpha_{\mathfrak{p}}\in\{0,12\}$ d'après la proposition~\ref{prop:ramification_caracteres} et on pose
\begin{equation*}
k=\sum_{\substack{\mathfrak{p}\mid p\\ \alpha_{\mathfrak{p}}=12 }}f_{\mathfrak{p}}\geq 0
\end{equation*}
de sorte que
\begin{equation}\label{eq:condition1}
\sum_{\mathfrak{p}\mid p}f_{\mathfrak{p}}\alpha_{\mathfrak{p}}=12k.
\end{equation}
Comme $p$ est non ramifié dans $K$, on a
\begin{equation}\label{eq:condition2}
d=\sum_{\mathfrak{p}\mid p}f_{\mathfrak{p}}.
\end{equation}
Or, d'après la remarque~\ref{rem:convention} suivant la proposition~\ref{prop:ramification_caracteres}, on peut toujours, si on le souhaite, remplacer la famille $\{\alpha_{\mathfrak{p}}\}_{\mathfrak{p}\mid p}$ par la famille $\{12-\alpha_{\mathfrak{p}}\}_{\mathfrak{p}\mid p}$, donc on peut supposer que l'on~a~:
\begin{equation*}
\sum_{\mathfrak{p}\mid p}f_{\mathfrak{p}}\alpha_{\mathfrak{p}}\leq \sum_{\mathfrak{p}\mid p}f_{\mathfrak{p}}(12-\alpha_{\mathfrak{p}}).
\end{equation*}
Autrement dit, d'après les égalités~(\ref{eq:condition1}) et~(\ref{eq:condition2})
\begin{equation*}
12k\leq 12(d-k)
\end{equation*}
soit encore
\begin{equation*}
k\leq \left[\frac{d}{2}\right].
\end{equation*}
D'après les égalités~(\ref{eq:RHS}) et~(\ref{eq:condition1}) on a
\begin{equation}\label{eq:RHS2}
\prod_{\mathfrak{p}\mid p}N_{\mathfrak{p}}(\ell+\mathfrak{p})^{\alpha_{\mathfrak{p}}}=\ell^{12k}\pmod{p}.
\end{equation}
Par ailleurs, d'après le lemme~\ref{lem:polynome_P_ell}, on a $\overline{P_{\ell}^*}(\Omega)=0$. Donc, d'après les égalités~(\ref{eq:demo_egalite_CFT}) et~(\ref{eq:RHS2}), il vient $\overline{P_{\ell}^*}\left(\ell^{12k}\right)=0\pmod{p}$, c'est-à-dire
\begin{equation*}
P_{\ell}^*\left(\ell^{12k}\right)\equiv 0\pmod{p}.
\end{equation*}
D'où le théorème~\ref{th:critere1}.

	\subsection{D\'emonstration du th\'eor\`eme~\ref{th:critere2}}

Soit $\ideal$ un id\'eal premier de bonne r\'eduction, $\gamma_{\ideal}$ un g\'en\'erateur de $\ideal^h$ et $\mathfrak{m}_{\gamma_{\ideal}}$ son polyn\^ome minimal sur~$\Qset$. On commence par un lemme pr\'eliminaire.
\begin{lemma}\label{lem:prelim_th2}
Le polyn\^ome $\mathfrak{m}_{\gamma_{\ideal}}^{(12)}$ appartient \`a $M_{\Zset}$ et v\'erifie pour $\mathfrak{p}$ divisant~$p$
\begin{equation*}
\overline{\left(\mathfrak{m}_{\gamma_{\ideal}}^{(12)}\right)^{*f_{\mathfrak{p}}}}\bigg(N_{\mathfrak{p}}(\gamma_{\ideal}+\mathfrak{p})^{12}\bigg)=0\pmod{p}.
\end{equation*}
\end{lemma}
\noindent{\it D\'emonstration. }Le polyn\^ome $\mathfrak{m}_{\gamma_{\ideal}}$ est irr\'eductible et unitaire. Il appartient donc \`a $M_{\Zset}$ et il en va de m\^eme pour $\mathfrak{m}_{\gamma_{\ideal}}^{(12)}$ d'apr\`es le lemme~\ref{lem:morphisme_monoide}. Par d\'efinition, on a 
\begin{equation*}
N_{\mathfrak{p}}(\gamma_{\ideal}+\mathfrak{p})=\gamma_{\ideal}\cdot\gamma_{\ideal}^p\cdots\gamma_{\ideal}^{p^{f_{\mathfrak{p}}-1}}+\mathfrak{p}\in\Zset/p\Zset
\end{equation*}
et
\begin{equation*}
\overline{\mathfrak{m}_{\gamma_{\ideal}}}(\gamma_{\ideal}+\mathfrak{p})=0\pmod{\mathfrak{p}}.
\end{equation*}
Or, le polyn\^ome $\mathfrak{m}_{\gamma_{\ideal}}$ est \`a coefficients dans~$\Zset$, d'o\`u $\overline{\mathfrak{m}_{\gamma_{\ideal} }^{(p)} }=\overline{\mathfrak{m}_{\gamma_{\ideal}} }$. On en d\'eduit donc avec les lemmes~\ref{lem:loi_monoide} et~\ref{lem:morphisme_monoide} que l'on a~:
\begin{equation*}
\begin{array}{clcl}
\overline{\mathfrak{m}_{\gamma_{\ideal} }^{(12)} }\big(\gamma_{\ideal}^{12}+\mathfrak{p}\big) & = & 0 \pmod{\mathfrak{p}} \\
\vdots & & \vdots \\
\overline{\mathfrak{m}_{\gamma_{\ideal} }^{(12)} }\big(\gamma_{\ideal}^{12p^{f_{\mathfrak{p}}-1}}+\mathfrak{p}\big) & = & 0 \pmod{\mathfrak{p}}.
\end{array}
\end{equation*}
Puis avec le lemme~\ref{lem:morphisme_anneaux_morphisme_monoides}, il vient
\begin{equation*}
\begin{array}{rcl}
& &  \overline{\left(\mathfrak{m}_{\gamma_{\ideal}}^{(12)}\right)^{*f_{\mathfrak{p}}}}\bigg(N_{\mathfrak{p}}(\gamma_{\ideal}+\mathfrak{p})^{12}\bigg) \\ 
& = &\overline{\mathfrak{m}_{\gamma_{\ideal}}^{(12)}}^{*f_{\mathfrak{p}}} \bigg(\gamma_{\ideal}^{12}\cdot\gamma_{\ideal}^{12p}\cdots\gamma_{\ideal}^{12p^{f_{\mathfrak{p}}-1}}+\mathfrak{p}\bigg) \\ 
& = & 0\pmod{p}
\end{array} 
\end{equation*}
car $N_{\mathfrak{p}}(\gamma_{\ideal}+\mathfrak{p})\in\Zset/p\Zset$. D'o\`u le lemme.

\bigskip

D\'emontrons \`a pr\'esent le th\'eor\`eme~\ref{th:critere2}. Supposons que $\ideal$ divise~$p$. Alors, $0$ est une racine commune de $P_{\ideal}^{(12h)}$ et $\mathfrak{m}_{\gamma_{\ideal}}^{(12)}$ modulo~$p$. Donc $p$ divise l'entier $\Res\left(P_{\ideal}^{(12h)},\mathfrak{m}_{\gamma_{\ideal}}^{(12)}\right)$ et par suite, si $d\geq 2$, $p$ divise $R_{\ideal}$.

Supposons que $\ideal$ ne divise pas~$p$. Alors, d'apr\`es la proposition~\ref{prop:CFT} appliqu\'ee \`a $a=\gamma_{\ideal}$, on a 
\begin{equation}\label{eq:egalite_CFT_th2}
\lambda(\sigma_{\ideal})^{12h}=\prod_{\mathfrak{p}\mid p}N_{\mathfrak{p}}\left(\gamma_{\ideal}+\mathfrak{p}\right)^{\alpha_{\mathfrak{p}}}= \prod_{\substack{\mathfrak{p}\mid p \\ \alpha_{\mathfrak{p}=12}}}N_{\mathfrak{p}}\left(\gamma_{\ideal}+\mathfrak{p}\right)^{12}
\end{equation}
Or, d'apr\`es le lemme~\ref{lem:prelim_th2}, on a
\begin{equation*}
\overline{\left(\mathfrak{m}_{\gamma_{\ideal}}^{(12)}\right)^{*f_{\mathfrak{p}}}}\bigg(N_{\mathfrak{p}}(\gamma_{\ideal}+\mathfrak{p})^{12}\bigg)=0\pmod{p}.
\end{equation*}
On en d\'eduit donc avec le lemme~\ref{lem:loi_monoide} que l'on a
\begin{equation*}
\myprod_{\substack{\mathfrak{p}\mid p \\ \alpha_{\mathfrak{p}=12}}} \overline{\left(\mathfrak{m}_{\gamma_{\ideal}}^{(12)}\right)^{*f_{\mathfrak{p}}}} \bigg(\prod_{\substack{\mathfrak{p}\mid p \\ \alpha_{\mathfrak{p}=12}}}N_{\mathfrak{p}}(\gamma_{\ideal}+\mathfrak{p})^{12}\bigg)=0\pmod{p},
\end{equation*}
puis avec l'\'egalit\'e~\eqref{eq:egalite_CFT_th2} ci-dessus et le lemme~\ref{lem:morphisme_anneaux_morphisme_monoides},
\begin{equation*}
\overline{\left(\mathfrak{m}_{\gamma_{\ideal}}^{(12)}\right)^{*k}} \left(\lambda(\sigma_{\ideal})^{12h}\right)=0\pmod{p},
\end{equation*}
o\`u l'on a pos\'e 
\begin{equation*}
k=\sum_{\substack{\mathfrak{p}\mid p\\ \alpha_{\mathfrak{p}}=12 }}f_{\mathfrak{p}}\geq 0.
\end{equation*}
Comme au \S\ref{sss:fin_demo_th1}, on peut supposer $k\leq[d/2]$. Par ailleurs, $\lambda(\sigma_{\ideal})^{12h}$ est une racine de $P_{\ideal}^{(12h)}\pmod{p}$. On en d\'eduit donc que $p$ divise $\Res\left(P_{\ideal}^{(12h)},\left(\mathfrak{m}_{\gamma_{\ideal}}^{(12)}\right)^{*k}\right)$ et par suite, $p$ divise~$R_{\ideal}$. Cela d\'emontre la premi\`ere partie du th\'eor\`eme~\ref{th:critere2}. Il reste \`a voir que si $E$ est sans multiplication complexe sur~$\overline{\Qset}$, alors $R_{\ideal}\not=0$ pour une infinit\'e de~$\ideal$. 

\medskip

Supposons $R_{\ideal}=0$. Alors, il existe une racine complexe $\alpha_{\ideal}$ de $P_{\ideal}$ telle que $\alpha_{\ideal}^{12h}$ soit racine de $\left(\mathfrak{m}_{\gamma_{\ideal}}^{(12)}\right)^{*k}$ pour un certain entier $0\leq k\leq [d/2]$. C'est impossible pour $k=0$ car $\alpha_{\ideal}^{12h}\not=1$. On a donc $k\geq 1$ (et par suite $d\geq 2$) et $\alpha_{\ideal}^{12h}$ s'\'ecrit comme un produit de~$k$ conjugu\'es de~$\gamma_{\ideal}$ \'elev\'es \`a la puissance~$12$. Notons $L^{\ideal}$ le corps engendr\'e par~$\alpha_{\ideal}$. C'est une extension de degr\'e au plus~$2$ de~$\Qset$. On distingue deux cas~:
\begin{enumerate}
	\item soit $\alpha_{\ideal}^{12h}\not\in\Qset$ et alors $L^{\ideal}=\Qset(\alpha_{\ideal}^{12h})$ est inclus dans $K^{\mathrm{gal}}$, la cl\^oture galoisienne de $K$ dans~$\overline{\Qset}$; en particulier, $L^{\ideal}$ est non ramifi\'e en dehors des premiers divisant~$D_K$.
	\item Soit $\alpha_{\ideal}^{12h}\in\Qset$ et alors 
\begin{equation*}
\zeta=\frac{\overline{\alpha_{\ideal}}}{\alpha_{\ideal}}
\end{equation*}
est une racine $12h$-i\`eme de l'unit\'e contenue dans~$L^{\ideal}$. C'est donc une racine primitive $2$-i\`eme, $3$-i\`eme, $4$-i\`eme ou $6$-i\`eme de l'unit\'e et l'on a~:
\begin{enumerate}
	\item soit $\zeta=1$, $t_{\ideal}^2=4N(\ideal)$ et $L^{\ideal}=\Qset$;
	\item soit $\zeta=-1$, $t_{\ideal}=0$ et $L^{\ideal}=\Qset(\sqrt{-1})$ ou $\Qset(\sqrt{-\ell})$;
	\item soit $\zeta=j$ ou $j^2$ (avec $j^2+j+1=0$), $t_{\ideal}^2=N(\ideal)$, donc $f_{\ideal}$ est pair et $L^{\ideal}=\Qset(\sqrt{-3})$;
	\item soit $\zeta=i$ ou $-i$ (avec $i^2=-1$), $t_{\ideal}^2=2N(\ideal)$, donc $\ell=2$ et $f_{\ideal}$ est impair. On en d\'eduit $L^{\ideal}=\Qset(\sqrt{-1})$;
	\item soit $\zeta=-j$ ou $-j^2$, $t_{\ideal}^2=3N(\ideal)$, donc $\ell=3$ et $f_{\ideal}$ est impair. On en d\'eduit et $L^{\ideal}=\Qset(\sqrt{-3})$.
\end{enumerate}
On en d\'eduit que dans ce cas la courbe $E$ a r\'eduction \emph{supersinguli\`ere} en $\ideal$ et que le corps $L^{\ideal}$ est non ramifi\'e en dehors de~$\{2,3,\ell\}$.
\end{enumerate}
Autrement dit, on a montr\'e que si $R_{\ideal}=0$, alors le corps $L^{\ideal}$ est non ramifi\'e en dehors des nombres premiers divisant $6\ell D_K$. Or, d'apr\`es un r\'esultat de Serre (\cite[IV-14(d)]{Ser68a}), on sait que si $E$ est sans multiplication complexe sur~$\overline{\Qset}$, alors pour tout ensemble fini $P$ de nombres premiers, il existe une infinit\'e d'id\'eaux premiers~$\ideal$ tels que $L^{\ideal}$ soit ramifi\'e en tout nombre premier 
appartenant \`a~$P$. Compte-tenu de l'\'etude pr\'ec\'edente, on en d\'eduit qu'il existe une infinit\'e d'id\'eaux premiers~$\ideal$ pour lesquels on a~$R_{\ideal}\not=0$. Cela ach\`eve la d\'emonstration du th\'eor\`eme~\ref{th:critere2}.

	\section{Bornes uniformes}

	\subsection{D\'emonstration de la proposition~\ref{prop:Phi_1}}

Soit $\mathfrak{q}$ un id\'eal premier de $\integers$ tel que le sous-groupe $\Phi_{\mathfrak{q}}$ soit non cyclique. Compte-tenu de la structure des groupes $\Phi$, l'id\'eal premier $\mathfrak{q}$ a nécessairement caractéristique résiduelle $\ell=2$ ou $3$ (\cite[\S$5.6(a)$]{Ser72}) et $\lvert\Phi_{\mathfrak{q}}\rvert=8$ ou $24$ (resp. $12$) si $\ell=2$ (resp. $\ell=3$). L'irréductibilité de $\rep$ résulte alors du lemme~\ref{lem:plgt_dans_Borel} ci-dessous et du fait que $\Phi_{\mathfrak{q}}$ se plonge dans $\mathrm{Aut}(\Etors)$ (car $\ell\not=p$ et $p\geq3$, \cite[\S$5.6(a)$]{Ser72}).

\medskip

\begin{notation}
On rappelle qu'un sous-groupe maximal de $\mathrm{Aut}(\Etors)$ stabilisant une droite de $\Etors$ est appel\'e sous-groupe de Borel.
\end{notation}

\begin{lemma}\label{lem:plgt_dans_Borel}
Soit $H$ un sous-groupe de $\mathrm{Gal}(K(\Etors)/K)$. On suppose $H$ non abélien fini. Si $p$ ne divise pas l'ordre de $H$, alors $H$ ne se plonge pas dans un sous-groupe de Borel de $\mathrm{Aut}(\Etors)$.
\end{lemma}
\noindent{\it D\'emonstration. }Supposons qu'il existe un morphisme injectif $\iota$ de $G$ dans un sous-groupe de Borel $B$ de $\mathrm{Aut}(\Etors)$. Dans une base convenable de $\Etors$ sur $\Fset_p$, $B$ est représentable matriciellement par le Borel standard
\begin{equation*}
\begin{pmatrix}
 * & * \\
 0  & * \\
\end{pmatrix}.
\end{equation*}
Il contient alors le sous-groupe $S$ d'ordre $p$ engendré par l'élément
\begin{equation*}
u=\begin{pmatrix}
 1 & 1 \\
 0  & 1 \\
\end{pmatrix}.
\end{equation*}
C'est un sous-groupe distingué de $B$. Comme l'ordre de $H$ est premier à $p$, le morphisme composé
\begin{equation*}
 H\stackrel{\iota}{\hookrightarrow} B\rightarrow B/S
\end{equation*}
est injectif. Par ailleurs, $B/S$ est abélien. D'où une contradiction car $H$ est supposé non abélien.

	\subsection{D\'emonstration de la proposition~\ref{prop:Phi_2}}

Soient $p\geq 3$ un nombre premier exceptionnel et $\mathfrak{q}$ un idéal premier de $\integers$ de caractéristique résiduelle $\ell\not=p$ en lequel $E$ a mauvaise réduction additive avec potentiellement bonne réduction. On souhaite montrer qu'il existe un entier $n\geq 0$ tel que l'ordre du groupe $\Phi_{\mathfrak{q}}$ divise $N(\mathfrak{q})^n(N(\mathfrak{q})-1)$. 

Vu la théorie du corps de classes, le caractère $\car$ s'interprète comme un homomorphisme 
\begin{equation*}
 \car: \mathrm{Gal}(K^{\mathfrak{m}}/K) \longrightarrow \Fset_p^*,
\end{equation*}
où $\mathfrak{m}$ est le conducteur de $\car$ et $K^{\mathfrak{m}}$ le corps de classes de rayon $\mathfrak{m}$. Alors, le caractère $\car$ est ramifié en $\mathfrak{q}$ (cf.~\cite[\S\S~1.12 et~5.6]{Ser72}) et on a une factorisation du type
\begin{equation*}
\mathfrak{m}=\mathfrak{m}'\cdot\mathfrak{q}^{n+1},\quad \textrm{où }n\geq 0\quad\textrm{ et }(\mathfrak{m}',\mathfrak{q})=1.
\end{equation*}
L'ordre du groupe $\Phi_{\mathfrak{q}}$ divise l'indice de ramification en $\mathfrak{q}$ de l'extension $K^{\mathfrak{m}}/K$. Or l'extension interm\'ediaire $K^{\mathfrak{m}'}/K$ est non ramifi\'ee en $\mathfrak{q}$. Donc l'ordre de $\Phi_{\mathfrak{q}}$ divise le cardinal du groupe $\mathrm{Gal}(K^{\mathfrak{m}}/K^{\mathfrak{m}'})$. Notons $h_{\mathfrak{m}}$ (resp. $h_{\mathfrak{m}'}$) le cardinal du groupe $ \mathrm{Gal}(K^{\mathfrak{m}}/K)$ (resp. $ \mathrm{Gal}(K^{\mathfrak{m}'}/K)$). Alors, d'après \cite[cor.$3.2.4$]{Coh00}, on a
\begin{equation*}
\lvert\mathrm{Gal}(K^{\mathfrak{m}}/K^{\mathfrak{m}'})\rvert=\frac{h_{\mathfrak{m}}}{h_{\mathfrak{m}'}} = \frac{(\mathcal{U}:\mathcal{U}_{\mathfrak{m}',1})}{(\mathcal{U}:\mathcal{U}_{\mathfrak{m},1})}N(\mathfrak{q})^n(N(\mathfrak{q})-1),
\end{equation*}
où $\mathcal{U}_{\mathfrak{m},1}$ (resp. $\mathcal{U}_{\mathfrak{m}',1}$) désigne le sous-groupe du groupe des unités $\mathcal{U}$ de $\integers$ qui sont congrues à $1$ modulo $\mathfrak{m}$ (resp. $\mathfrak{m}'$) au sens de \cite[Def.$3.2.2$]{Coh00}. Or, comme $\mathfrak{m}'$ divise $\mathfrak{m}$, l'indice de $\mathcal{U}_{\mathfrak{m}',1}$ dans $\mathcal{U}$ divise celui de $\mathcal{U}_{\mathfrak{m},1}$. Donc, l'ordre de $\mathrm{Gal}(K^{\mathfrak{m}}/K^{\mathfrak{m}'})$ divise $N(\mathfrak{q})^n(N(\mathfrak{q})-1)$ et il en va de m\^eme en particulier pour l'ordre de $\Phi_{\mathfrak{q}}$. D'où la proposition~\ref{prop:Phi_2}.

\medskip

\noindent{\it Remarque. }Lorsque $\ell\geq 5$, on a $\lvert\Phi_{\mathfrak{q}}\rvert=2$, $3$, $4$ ou $6$ (\cite[p. 312]{Ser72}). Or $N(\mathfrak{q})$ est premier à $12$, donc $\lvert\Phi_{\mathfrak{q}}\rvert$ divise $N(\mathfrak{q})^n(N(\mathfrak{q})-1)$ pour un certain entier $n$ \ssi $\lvert\Phi_{\mathfrak{q}}\rvert$ divise $N(\mathfrak{q})-1$. 
Cela justifie la remarque~\ref{rem:precisionpaumoins5}.

	\subsection{D\'emonstration des corollaires~\ref{cor:Phi_2} et~\ref{cor:Phi_3}}

Supposons que $\mathfrak{q}$ divise $2$. Lorsque $\lvert\Phi_{\mathfrak{q}}\rvert=8$ ou $24$, le groupe $\Phi_{\mathfrak{q}}$ n'est pas abélien (\cite[$5.6(a)$]{Ser72}) et la conclusion résulte  de la prop.~\ref{prop:Phi_1}. Pour $\lvert\Phi_{\mathfrak{q}}\rvert=3$ ou $6$, supposons la représentation $\rep$ réductible. Alors, d'après la prop.~\ref{prop:Phi_2}, l'ordre de $\Phi_{\mathfrak{q}}$ divise $2^{f_{\mathfrak{q}}}(2^{f_{\mathfrak{q}}}-1)$. Or, $2^{f_{\mathfrak{q}}}-1\equiv 1\pmod{3}$ car $f_{{\mathfrak{q}}}$ est impair. D'où une contradiction et le corollaire~\ref{cor:Phi_2}.

Supposons que $\mathfrak{q}$ divise $3$. Lorsque $\lvert\Phi_{\mathfrak{q}}\rvert=12$, le groupe $\Phi_{\mathfrak{q}}$ n'est pas abélien (\cite[$5.6(a)$]{Ser72}) et la conclusion résulte comme ci-dessus de la prop.~\ref{prop:Phi_1}. Pour $\lvert\Phi_{\mathfrak{q}}\rvert=4$, supposons la représentation $\rep$ réductible. Alors, d'après la prop.~\ref{prop:Phi_2}, l'ordre de $\Phi_{\mathfrak{q}}$ divise $3^{f_{\mathfrak{q}}}(3^{f_{\mathfrak{q}}}-1)$. Or, $3^{f_{\mathfrak{q}}}-1\equiv 2\pmod{4}$ car $f_{{\mathfrak{q}}}$ est impair, d'où une contradiction et le corollaire~\ref{cor:Phi_3}.

	\section{Exemples numériques}\label{s:exemples}

L'objet de cette section est de d\'eterminer explicitement, pour certaines courbes elliptiques~$E$ définies sur des corps de nombres $K$, l'ensemble $\Exc[E/K]$ et d'illustrer ainsi chacun des r\'esultats du \S\ref{s:enonces_resultats}.

	\subsection{Stratégie}
On traite successivement des courbes d\'efinies sur des corps quadratiques, puis une sur un corps cubique et enfin une derni\`ere sur un corps biquadratique. Pour la plupart d'entre elles, on commence par déterminer son type de réduction en chaque idéal premier. Si pour l'un d'entre eux, on est dans un cas d'application des \og résultats uniformes\fg\ (prop.~\ref{prop:Phi_1} et~\ref{prop:Phi_2} et cor.~\ref{cor:Phi_2} et~\ref{cor:Phi_3}) du~\S\ref{s:enonces_resultats}, il ne reste plus alors à traiter que le cas $p=2$ et éventuellement $p=3$ et $p=\ell$ où $\ell$ est un nombre premier $\geq 5$. Sinon, on applique le critère du théorème~\ref{th:critere1}. On cherche alors un nombre premier $\ell$ de bonne réduction pour lequel $B_{\ell}$ soit non nul. Si $d$ est impair, c'est automatique (cor.~\ref{cor:degre_impair}). Si $d$ est pair, l'existence d'un tel nombre premier n'est malheureusement pas garantie comme le montre l'exemple~\ref{ex:biquad}. Cependant dans le cas quadratique ($d=2$), cela ne pose aucun problème d'en trouver un dans la pratique. 

Après quelques itérations, on obtient alors un ensemble très restreint de nombres premiers contenant $2$, $3$ et les premiers ramifiés dans le corps. On traite ensuite \og à la main\fg\ ceux qui restent. Soit on trouve un idéal premier $\mathfrak{q}$ de bonne réduction ne divisant pas~$p$ tel que $P_{\mathfrak{q}}$ soit irréductible modulo~$p$ et alors $p$ n'est pas exceptionnel, soit on montre que $E$ possède un sous-groupe stable d'ordre~$p$ et alors $p$ est exceptionnel.

		\subsection{Notations}

La courbe $E$ est donn\'ee sous forme d'une \'equation de Weierstrass
\begin{equation*}
y^2+a_1xy+a_3y=x^3+a_2x^2+a_4x+a_6,
\end{equation*}
avec $a_i\in\integers$. On adopte les notations standard de Tate (\cite{Tat75}). Pour chaque id\'eal premier $\mathfrak{p}$ de $\integers$, on note $v_{\mathfrak{p}}$ la valuation en $\mathfrak{p}$ de $K$ normalis\'ee par $v_{\mathfrak{p}}(K^*)=\Zset$.

\'Etant donné un nombre premier~$\ell$, on note $\mathfrak{S}_{\ell}^{(d)}$ l'ensemble des 
diviseurs premiers de $B_{\ell}$. Lorsque $E$ est d\'efinie sur un corps quadratique, on dispose 
de deux programmes \verb"pari", disponibles à l'adresse~: 
\begin{center}
\verb"http://people.math.jussieu.fr/~billerey/".
\end{center}
appel\'es, \verb"TraceOfFrobenius" et \verb"ExceptionalPrimes", permettant respectivement de calculer la famille $\{t_{\ideal}\}_{\ideal\mid\ell}$ et de déterminer l'ensemble $\mathfrak{S}_{\ell}^{(2)}$ (et aussi l'entier $B_{\ell}$).

		\subsection{Les polyn\^omes $P_{\ell}^*$ dans le cas quadratique}\label{ss:cas_quad}

On suppose que $K$ est un corps quadratique, i.e. $d=2$. Pour un $\ell$ de bonne r\'eduction, on a alors $B_{\ell}=P_{\ell}^*(1)\cdot P_{\ell}^*(\ell^{12})$ et $P_{\ell}^*(1)\not=0$ (lemme~\ref{lem:polynome_P_ell}). La proposition suivante donne une description explicite des polyn\^omes $P_{\ell}^*$ et de la condition $B_{\ell}=0$. On rappelle au préalable que pour tout entier $n\geq 1$, il existe un unique polynôme $T_n$ appartenant à $\Zset[X]$ tel que pour tout nombre réel $\theta$, on ait $T_n(\cos\theta)=\cos(n\theta)$. Le polynôme $T_n$ s'appelle le $n$-ième polynôme de Tchebychev (de première espèce). On a en particulier,
\begin{equation*}
T_{12}(X)=2048X^{12} - 6144X^{10} + 6912X^8 - 3584X^6 + 840X^4 - 72X^2 + 1
\end{equation*}
et $T_{24}(X)=2T_{12}(X)^2-1$.

\begin{proposition}\label{prop:polynomes_explicites}
On suppose $K/\Qset$ quadratique. Soit $\ell$ nombre premier de bon\-ne réduction. On est dans l'une des situations suivantes.
\begin{enumerate}
	\item Soit $\ell$ est ramifié dans $K$, $\ell\integers=\ideal^2$ et on a
\begin{equation*}
P_{\ell}^*(X) = P_{\ideal}^{(24)}(X)=X^2-2\ell^{12}T_{24}\left(t_{\ideal}/2\sqrt{\ell}\right)X+\ell^{24}.
\end{equation*}
En particulier, on a
\begin{equation*}
P_{\ell}^*(\ell^{12})=-\ell^{12}t_{\ideal}^2(t_{\ideal}^2-\ell)^2(t_{\ideal}^2-4\ell)(t_{\ideal}^2-2\ell)^2(t_{\ideal}^2-3\ell)^2(t_{\ideal}^4-4\ell t_{\ideal}^2+\ell^2)^2.
\end{equation*}
Ainsi $P_{\ell}^*(\ell^{12})=0$ \ssi $t_{\ideal}\equiv 0\pmod{\ell}$, c'est-à-dire \ssi $E$ a bonne r\'eduction \emph{supersinguli\`ere} en~$\ideal$.
	\item Soit $\ell$ est inerte dans $K$, $\ell\integers=\ideal$ et on a
\begin{equation*}
P_{\ell}^*(X)= P_{\ideal}^{(12)}(X)=X^2-2\ell^{12}T_{12}\left(t_{\ideal}/2\ell\right)X+\ell^{24}.
\end{equation*}
En particulier, on a
\begin{equation*}
P_{\ell}^*(\ell^{12})=-\ell^{12}t_{\ideal}^2(t_{\ideal}^2-\ell^2)^2(t_{\ideal}^2-4\ell^2)(t_{\ideal}^2-3\ell^2)^2.
\end{equation*}
Ainsi $P_{\ell}^*(\ell^{12})=0$ \ssi $t_{\ideal}\equiv 0\pmod{\ell}$, c'est-à-dire \ssi $E$ a bonne r\'eduction \emph{supersinguli\`ere} en~$\ideal$.
	\item Soit $\ell$ est décomposé dans $K$, $\ell\integers=\ideal_1\ideal_2$ et on a
\begin{multline*}
 P_{\ell}^*(X)= (P_{\ideal_1}\ast P_{\ideal_2})^{(12)}(X)
= X^4-4\ell^{12}T_{12}(t_{\ideal_1}/2\sqrt{\ell})T_{12}(t_{\ideal_2}/2\sqrt{\ell})X^3 \\
 -2\ell^{24}(1-2\left(T_{12}\left(t_{\ideal_1}/2\sqrt{\ell})^2+T_{12}(t_{\ideal_2}/2\sqrt{\ell})^2\right)\right)X^2 \\
 -4\ell^{36}T_{12}(t_{\ideal_1}/2\sqrt{\ell})T_{12}(t_{\ideal_2}/2\sqrt{\ell})X+\ell^{48}.
\end{multline*}
En particulier, on a
\begin{multline*}
P_{\ell}^*(\ell^{12})=\ell^{36}\left(t_{\ideal_1}^2-t_{\ideal_2}^2\right)^2\left((t_{\ideal_1}^2+t_{\ideal_2}^2-3\ell)^2-t_{\ideal_1}^2t_{\ideal_2}^2\right)^2\left(t_{\ideal_1}^2+t_{\ideal_2}^2-4\ell\right)^2\\
\times\left((t_{\ideal_1}^2+t_{\ideal_2}^2-\ell)^2-3t_{\ideal_1}^2t_{\ideal_2}^2\right)^2.
\end{multline*}
Ainsi $P_{\ell}^*(\ell^{12})=0$ \ssi l'une des conditions suivantes est satisfaite~:
\begin{equation*}
t_{\ideal_1}=\pm t_{\ideal_2};\quad t_{\ideal_1}^2+t_{\ideal_2}^2\pm t_{\ideal_1}t_{\ideal_2}=3\ell;\quad t_{\ideal_1}^2+t_{\ideal_2}^2=4\ell.
\end{equation*}
\end{enumerate}
\end{proposition}
\noindent{\it D\'emonstration. }La preuve de cette proposition repose sur la proposition~\ref{prop:RH} ainsi que sur les relations de r\'ecurrence entre polyn\^omes de Tchebychev. Elle n'est pas difficile. On ne traite que le cas o\`u $\ell$ est inerte, les autres \'etant analogues. Supposons donc $\ell$ inerte dans $K$ avec $\ell\integers=\ideal$ et posons
\begin{equation*}
P_{\ideal}(X)=X^2-t_{\ideal}X+\ell^2=(X-\alpha)(X-\beta).
\end{equation*}
D'après la prop.~\ref{prop:RH}, on a $\lvert\alpha\rvert=\lvert\beta\rvert=\ell$. Posons donc $\alpha=\ell e^{i\theta}$ avec $\theta\in\Rset$.
D'après le lemme~\ref{lem:loi_monoide}, on a
\begin{equation*}
P_{\ell}^*(X)=(X-\alpha^{12})(X-\beta^{12}).
\end{equation*}
D'où
\begin{equation*}
P_{\ell}^*(X)= X^2-(\alpha^{12}+\beta^{12})X+\ell^{24} = X^2-2\ell^{12}\cos(12\theta)X+\ell^{24}.
\end{equation*}
Or, $\cos(12\theta)=T_{12}(\cos\theta)$ et $2\ell\cos\theta=t_{\ideal}$, d'où
\begin{equation*}
P_{\ell}^*(X) = X^2-2\ell^{12}T_{12}\left(\frac{t_{\ideal}}{2\ell}\right)X+\ell^{24}.
\end{equation*}
On en déduit immédiatement
\begin{equation*}
P_{\ell}^*(\ell^{12})= 2\ell^{24}\left(1-T_{12}\left(\frac{t_{\ideal}}{2\ell}\right)\right).
\end{equation*}
Or, on a $1-T_{12}=8(1-T_3)(1+T_3)T_3^2$. D'où la factorisation
\begin{equation*}
P_{\ell}^*(\ell^{12})=-\ell^{12}t_{\ideal}^2(t_{\ideal}^2-\ell^2)^2(t_{\ideal}^2-4\ell^2)(t_{\ideal}^2-3\ell^2)^2
\end{equation*}
car 
\begin{equation*}
T_3(X)=4X^3-3X;\quad 1-T_3(X)=-(X-1)(2X+1)^2;
\end{equation*}
et
\begin{equation*}
1+T_3(X)=(X+1)(2X-1)^2.
\end{equation*}
On en déduit que l'on a $P_{\ell}^*(\ell^{12})=0$ \ssi $t_{\ideal}=0$, $\pm\ell$ ou $\pm 2\ell$. Autrement dit, \ssi $t_{\ideal}\equiv 0\pmod{\ell}$ car $\lvert t_{\ideal}\rvert\leq 2\ell$.

	\subsection{Exemples}
On illustre dans les exemples suivants \emph{chacun} des r\'esultats du \S\ref{s:enonces_resultats}.

\begin{example}
On suppose $K=\Qset(\sqrt{5})$. On consid\`ere la courbe $E$ d'\'equation
\begin{equation*}
y^2=x^3+2x^2+\omega x\quad\textrm{où}\quad \omega=\frac{1+\sqrt{5}}{2}.
\end{equation*}
Alors, $\Exc[E/K]=\{2\}$.
\end{example}
\noindent{\it D\'emonstration. }On a
\begin{equation*}
\left\{
\begin{array}{l}
c_4=2^4(4-3\omega) \\
c_6=2^6(-8+9\omega) \\
\Delta=-2^6\omega. \\
\end{array}
\right.
\end{equation*}
Or, $\omega$ est une unit\'e de $\integers$. En particulier, la courbe $E$ a bonne r\'eduction en dehors de (l'id\'eal premier) $2\integers$. On a~:
\begin{equation*}
(v_{2}(c_4),v_{2}(c_6),v_{2}(\Delta))=(4,6,6).
\end{equation*}
Donc $E$ a mauvaise r\'eduction additive en~$2$. On note $\Phi_2$ son défaut de semi-stabilité en~$2$. On a $v_{2}(j)=6$ et $3v_{2}(c_4)=2v_{2}(c_6)$. L'extension $K/\Qset$ \'etant non ramifié en~$2$, on a d'apr\`es~\cite{Cal04}, $\lvert\Phi_{2}\rvert=4$ ou~$8$. Or, avec les notations de loc.~cit. la condition~(C2) n'est pas satisfaite. On en d\'eduit que l'on a $\lvert\Phi_{2}\rvert=8$. Et, d'apr\`es le cor.~\ref{cor:Phi_2}, $\rep$ est irr\'eductible pour tout nombre premier $p\geq 5$. La courbe $E$ a bonne réduction en l'idéal premier $7\integers$ et d'apr\`es le programme \verb"TraceOfFrobenius", on a $t_7=-12$. D'où
\begin{equation*}
P_{7}(X)=X^2-t_7X+49\equiv X^2+1\pmod{3}.
\end{equation*}
Donc $\rep[3]$ est \'egalement irr\'eductible. La repr\'esentation $\rep[2]$, en revanche, est r\'eductible car $(0,0)$ est un point d'ordre~$2$.

\begin{example}
On suppose $K=\Qset(\sqrt{13})$. On consid\`ere la courbe $E$ d'\'equation
\begin{equation*}
y^2=x^3-(313+240\omega)x-17\quad\textrm{où}\quad\omega=\frac{1+\sqrt{13}}{2}.
\end{equation*}
Alors, l'ensemble $\Exc[E/K]$ est vide.
\end{example}
\noindent{\it D\'emonstration. }On a
\begin{equation*}
\left\{
\begin{array}{l}
c_4=2^4\cdot 3(11+8\omega)^2 \\
c_6=2^5\cdot 3^3\cdot 17 \\
\Delta=2^4\cdot 5(11+8\omega)^2(213629+167568\omega). \\
\end{array}
\right.
\end{equation*}
De plus, $N_{K/\Qset}(213629+167568\omega)=-1153\cdot 2430503$ et ni $1153$, ni $2430503$ ne divisent $c_4$. Donc la courbe $E$ a mauvaise réduction multiplicative en un idéal premier au-dessus de $1153$ et un idéal premier au-dessus de $2430503$. Le nombre premier $2$ est inerte dans $K$ et 
\begin{equation*}
(v_2(c_4),v_2(c_6),v_2(\Delta))=(4,5,4).
\end{equation*}
Donc $v_2(j)=8$ et d'après~\cite{Cal04}, on a $\lvert\Phi_2\rvert=3$, $6$ ou $24$. Comme par ailleurs, 
\begin{equation*}
j'=\frac{j}{2^8}\equiv -1\pmod{4},
\end{equation*}
la condition~(C3) de \emph{loc.~cit. }est satisfaite avec $\gamma=1$ et $\lvert\Phi_2\rvert=3$ ou~$6$ (en fait $\lvert\Phi_2\rvert=6$ d'après loc.~cit.). Puisque $f_2=2$ est pair, le cor.~\ref{cor:Phi_2} ne s'applique pas. Cependant, en l'idéal premier $\ideal_{17}=(15+4\sqrt{13})\integers$, on a 
\begin{equation*}
(v_{\ideal_{17}}(c_4),v_{\ideal_{17}}(c_6),v_{\ideal_{17}}(\Delta))=(2,1,2).
\end{equation*}
Donc $E$ a mauvaise réduction additive en $\ideal_{17}$ avec potentiellement bonne réduction. Son défaut de semi-stabilité $\Phi_{\ideal_{17}}$ est d'ordre~$6$ (\cite[p.312]{Ser72}). Or, $6$ ne divise pas $N(\ideal_{17})-1=16$. Donc, d'après la prop.~\ref{prop:Phi_2}, la représentation~$\rep$ est irréductible pour tout nombre premier $p\geq 3$ et $p\not=17$. Si $\ideal_3$ désigne un idéal divisant $3$, alors $E$ a bonne réduction en~$\ideal_3$ et d'après le programme \verb"TraceOfFrobenius", on a
\begin{equation*}
t_{\ideal_3}=-3.
\end{equation*}
Donc le polynôme $P_{\ideal_3}(X)=X^2+3X+3$ est irréductible modulo~$2$ et~$17$. On en déduit le résultat.

\begin{example}
On suppose $K=\Qset(\sqrt{-1})$. Pour tout entier $a\in\Zset[\sqrt{-1}]$, on consid\`ere la courbe $E_a$ d'\'equation
\begin{equation*}
y^2=x^3+ax+a.
\end{equation*}
Supposons que l'on ait $v_{\ideal_2}(a)=0,1,3,4,5,6$ ou $7$, où $\ideal_2$ est l'unique idéal de $\integers$ au-dessus de~$2$. Alors, l'ensemble $\Exc[E_a/K]$ est contenu dans $\{2,3\}$.
\end{example}
\noindent{\it D\'emonstration. }On a
\begin{equation*}
\left\{
\begin{array}{l}
c_4(E_a)=-2^4\cdot 3\cdot a \\
c_6(E_a)=-2^5\cdot 3\cdot a \\
\Delta(E_a)=-2^4\cdot a^2(4a+27). \\
\end{array}
\right.
\end{equation*}
En particulier, $v_{\ideal_2}(c_4(E_a))=8+v_{\ideal_2}(a)$, $v_{\ideal_2}(\Delta(E_a))=8+2v_{\ideal_2}(a)$ et
\begin{equation*}
v_{\ideal_2}(j(E_a))=16+v_{\ideal_2}(a)\geq 0. 
\end{equation*}
Donc $E_a$ a mauvaise réduction additive en $\ideal_2$ avec potentiellement bonne réduction. Par ailleurs, $2$ est ramifié dans $K$ (donc $f_{\ideal_2}=1$ est impair) et
\begin{equation*}
v_{\ideal_2}(j(E_a))\in\{16,17,19,20,21,22,23\}. 
\end{equation*}
Donc, d'après \cite{Bil09a}, on a $\lvert\Phi_{\ideal_2}\rvert\in\{3,6,8,24\}$. On conclut avec le cor.~\ref{cor:Phi_2}.

\medskip

Avec les notations de l'exemple précédent, lorsque $v_{\ideal_2}(a)=2$ ou $\geq 8$, le corollaire~\ref{cor:Phi_2} ne s'applique pas toujours. On traite ci-dessous un exemple dans le cas où $v_{\ideal_2}(a)=2$.

\begin{example}
On suppose $K=\Qset(\sqrt{-1})$. On consid\`ere la courbe $E$ d'\'equation
\begin{equation}\label{eq:exemple_sur_Q(i)}
y^2=x^3+2(3+2\sqrt{-1})x+2(3+2\sqrt{-1}).
\end{equation}
Alors, l'ensemble $\Exc[E/K]$ est vide.
\end{example}
\noindent{\it D\'emonstration. }On a
\begin{equation*}
\left\{
\begin{array}{l}
c_4=-2^5\cdot 3(3+2\sqrt{-1}) \\
c_6=-2^6\cdot 3^3(3+2\sqrt{-1}) \\
\Delta=-2^6(3+2\sqrt{-1})^2\cdot (51+16\sqrt{-1}). \\
\end{array}
\right.
\end{equation*}
On a $2\integers=\ideal_2^2$ o\`u $\ideal_2=(1+\sqrt{-1})\integers$ et
\begin{equation*}
(v_{\ideal_2}(c_4),v_{\ideal_2}(c_6),v_{\ideal_2}(\Delta))=(10,12,12).
\end{equation*}
Donc, d'apr\`es \cite{Pap93}, l'\'equation~(\ref{eq:exemple_sur_Q(i)}) correspond \`a un cas $6$ ou $7$ de Tate. En particulier, il est minimal en $\ideal_2$ et $E$ a r\'eduction additive en $\ideal_2$. D\'eterminons \`a pr\'esent l'ordre $\lvert\Phi_{\ideal_2}\rvert$ de son d\'efaut de semi-stabilit\'e  en $\ideal_2$. On a $v_{\ideal_2}(j)=18$ et $\pi=1+\sqrt{-1}$ est une uniformisante de $K_{\ideal_2}$. De plus,
\begin{equation*}
\frac{c_4}{\pi^{10}}\equiv 1+\pi\pmod{2}.
\end{equation*}
Donc, d'apr\`es \cite[th.~2]{Bil09a}, on a $\lvert\Phi_{\ideal_2}\rvert=4$. On ne peut donc pas appliquer le cor.~\ref{cor:Phi_2}.

Notons $\ideal_{13}$ l'id\'eal premier de $\integers$ engendr\'e par $3+2\sqrt{-1}$. On a
\begin{equation*}
(v_{\ideal_{13}}(c_4),v_{\ideal_{13}}(\Delta))=(1,2),
\end{equation*}
d'où $v_{\ideal_{13}}(j)=1$. L'\'equation~(\ref{eq:exemple_sur_Q(i)}) est minimale en $\ideal_{13}$ et $E$ a mauvaise r\'eduction additive en $\ideal_{13}$ avec potentiellement bonne réduction. Son d\'efaut de semi-stabilit\'e est d'ordre $6$ (c.f.~\cite[p.312]{Ser72}). Comme $6$ divise $N(\ideal_{13})-1=12$, la prop.~\ref{prop:Phi_2} ne donne aucune majoration de $\Exc[E/K]$. Par ailleurs, on a $v_{\ideal_{13}}(\Delta)=2$ et d'après \cite[lem.1]{Kra97a}, $E$ n'a pas potentiellement bonne réduction de hauteur~$2$ en~$\ideal_{13}$.

Notons $\ideal_{2857}$ l'id\'eal premier de $\integers$ engendr\'e par $51+16\sqrt{-1}$. La courbe $E$ a mauvaise r\'eduction multiplicative en $\ideal_{2857}$. En dehors de  $\ideal_{2}$,  $\ideal_{13}$ et  $\ideal_{2857}$, la courbe $E$ a bonne r\'eduction.

Vu l'\'etude pr\'ec\'edente, aucun des r\'esultats \og uniformes\fg\ du \S\ref{s:enonces_resultats} ne s'applique. Pour cette courbe, on a donc recours au crit\`ere du th.~\ref{th:critere1}. Soit $p\geq 5$ un nombre premier exceptionnel. D'apr\`es le programme \verb"TraceOfFrobenius", on a
\begin{equation*}
\{t_{\ideal}\}_{\ideal\mid 5}=\{-2,1\}\quad\textrm{et}\quad t_7=6.
\end{equation*}
D'apr\`es le th.~\ref{th:critere1} appliqu\'ee \`a $\ell=5$ et $\ell=7$, $p$ divise chacun des entiers $B_{5}$ et $B_{7}$. Or, d'apr\`es le programme \verb"ExceptionalPrimes", on a
\begin{equation*}
B_{5}=2^{28}\cdot 3^{16}\cdot 5^{39}\cdot 11^2\cdot 17\cdot 61\cdot 73\cdot 277\cdot 397\cdot 557\cdot 653\cdot 757\cdot 23833
\end{equation*}
et
\begin{equation*}
B_{7}=2^{14}\cdot 3^8\cdot 5^2\cdot 7^{13}\cdot 11\cdot 13^5\cdot 37^2\cdot 2089\cdot 2689\cdot 3889.
\end{equation*}
D'où
\begin{equation*}
p\in \mathfrak{S}_{5}^{(2)}\cap \mathfrak{S}_{7}^{(2)}=\{2,3,5,11\}.
\end{equation*}
Il ne reste donc plus qu'\`a traiter les cas $p=2$, $3$, $5$ et $11$. Or, $E$ a bonne réduction en l'idéal premier $3\integers$ et d'après le programme \verb"TraceOfFrobenius", on a
\begin{equation*}
P_{3}(X)=X^2+3X+9.
\end{equation*}
Donc $P_{3}$ est irr\'eductible modulo~$2$, $5$ et $11$. Et, si $\ideal_{5}$ est un id\'eal premier au-dessus de~$5$, on a $t_{\ideal_{5}}=-2$ ou $1$, et
\begin{equation*}
P_{\ideal_{5}}(X)\equiv X^2+2X+2\pmod{3}.
\end{equation*}
Donc $P_{\ideal_{5}}$ est irr\'eductible modulo~$3$. On en d\'eduit que $\rep$ est \'egalement irr\'eductible pour $p=2$, $3$, $5$ et $11$. D'o\`u le r\'esultat.

\begin{example}
On suppose $K=\Qset(\sqrt{2})$ et on pose
\begin{equation*}
\left\{\begin{array}{lcl}
	A & = & -3^3\cdot 5\cdot 17^3(428525+303032\sqrt{2}) \\
	B & = & 2\cdot 3^3\cdot 5\cdot 17^3(62176502533+43965551956\sqrt{2}).
\end{array}
\right.
\end{equation*}
On consid\`ere la courbe $E$ d'\'equation
\begin{equation*}
y^2=x^3+Ax+B.
\end{equation*}
Alors, $\Exc[E/K]=\{13\}$.
\end{example}
\noindent{\it D\'emonstration. }On vérifie que pour le modèle choisi, on a
\begin{equation*}
N_{K/\Qset}(\Delta)=-2^{25}\cdot 3^{18}\cdot 5^{4}\cdot 7^{2}\cdot 17^{15}\cdot 23^{6}\cdot 79^{6}.
\end{equation*}
En particulier, la courbe $E$ a bonne réduction en les idéaux premiers divisant $11$, $13$, $19$, $29$ et $41$ et d'après le programme \verb"TraceOfFrobenius", on a
\begin{equation*}
t_{11}=4;\quad t_{13}=-14\quad t_{19}=26; \quad t_{29}=1\quad\textrm{et}\quad \{t_{\ideal}\}_{\ideal\mid 41}=\{-3,2\}.
\end{equation*}
Soit $p$ un nombre premier exceptionnel n'appartenant pas à l'ensemble
\begin{equation*}
\{2,3,5,7,17,23,79\}.
\end{equation*}
Alors, d'après le corollaire~\ref{cor:discriminant}, on a en particulier
\begin{equation*}
p\in \mathfrak{S}_{11}^{(2)}\cap \mathfrak{S}_{13}^{(2)}.
\end{equation*}
D'où, d'après le programme \verb"ExceptionalPrimes", 
\begin{equation*}
p\in \{2,3,5,7,13\}.
\end{equation*}
Autrement dit, il ne reste plus qu'à traiter les cas où $p=2,3,5,7,13,17,23$ et $79$. Or le polynôme $P_{11}$ est irréductible modulo~$5$, $23$ et $79$. De même, $P_{13}$ est irréductible modulo~$7$, $P_{19}$ modulo~$17$ et $P_{29}$ modulo~$2$. Si $\ideal_{41}$ désigne l'idéal premier de $\integers$ au-dessus de~$41$ tel que $t_{\ideal_{41}}=2$, alors $P_{\ideal_{41}}$ est irréductible modulo~$3$. On en déduit que $2,3,5,7,17,23$ et $79$ ne sont pas exceptionnels. En revanche $13$ est un nombre premier exceptionnel. En effet, la courbe modulaire $X_0(13)$ paramétrisant les courbes elliptiques munies d'un sous-groupe stable d'ordre~$13$ est de genre~$0$ et un isomorphisme avec $\mathbf{P}^1$ est donné par la fraction rationnelle suivante (\cite[\S2.2]{Mes80})~:
\begin{equation*}
j(X_0(13))(X)=\frac{(X^2+5X+13)(X^4+7X^3+20X^2+19X+1)^3}{X}.
\end{equation*}
On vérifie alors que l'on a $j(X_0(13))(\sqrt{2})=j$. Cela montre que $13$ est exceptionnel et le résultat annoncé.

\begin{example}
On suppose $K=\Qset(\sqrt{3})$ et on pose
\begin{equation*}
\left\{\begin{array}{lcl}
	a_1 & = & 2^2\cdot 7\sqrt{3}(1+2\sqrt{3})(2-\sqrt{3})=252-112\sqrt{3} \\
	a_4 & = & 2\cdot 3^2\cdot 7^2\sqrt{3}(1+2\sqrt{3})^2=10584+11466\sqrt{3} \\
	a_6 & = & 2^3\cdot 3\cdot 7^4\sqrt{3}(1+2\sqrt{3})^4(7-4\sqrt{3})=- 24202080+15616104\sqrt{3}. \\
\end{array}
\right.
\end{equation*}
On consid\`ere la courbe $E$ d'\'equation
\begin{equation}\label{ex:courbeCM_avec_Exc_fini}
y^2+a_1xy=x^3+a_4x+a_6.
\end{equation}
Alors, la courbe $E$ a des multiplications complexes par le corps $\Qset(\sqrt{-1})$ et $\Exc[E/K]=\{2,3\}$.
\end{example}
\noindent{\it D\'emonstration. }On a
\begin{equation*}
\left\{
\begin{array}{rl}
c_4 & =- 2^5\cdot 3^2\cdot 7^2\cdot 11^2(-1+2\sqrt{3})(2-3\sqrt{3})^3\varepsilon^{-2} \\
c_6 & = 2^9\cdot 3^3\cdot 7^4\cdot 11^3(1+2\sqrt{3})(2-3\sqrt{3})^3\varepsilon^{-4} \\
\Delta & =- 2^9\cdot 3^4\cdot 7^6\cdot 11^6\sqrt{3}(2-3\sqrt{3})^6\varepsilon^{-4} \\
\end{array}
\right.
\end{equation*}
où $\varepsilon=2+\sqrt{3}$ est l'unité fondamentale de $\Qset(\sqrt{3})$. On en déduit que
\begin{align*}
j & =2^6\cdot 3\cdot \sqrt{3}(-1+2\sqrt{3})^3(2-3\sqrt{3})^3\varepsilon^{-2} \\
  & = 76771008-44330496\sqrt{3}
\end{align*}
est entier de polyn\^ome minimal sur $\Qset$
\begin{equation*}
P(X)=X^2-153542016X-1790957481984.
\end{equation*}
Vérifions que $E$ a des multiplications complexes par l'ordre de $\Qset(\sqrt{-1})$ de conduteur~$3$. En effet, il n'y a qu'un nombre fini de classes d'isomorphisme de courbes elliptiques ayant des multiplications complexes par un ordre de discriminant $D$ fix\'e. De plus, le polyn\^ome minimal sur $\Qset$ de l'invariant modulaire d'une telle courbe elliptique est
\begin{equation*}
\Phi_D(X)=\prod_{a,b,c}\left(X-j\left(\frac{-b+\sqrt{D}}{2a}\right)\right),
\end{equation*}
où $j$ d\'esigne la fonction modulaire et $(a,b,c)$ parcourt l'ensemble des triplets d'entiers tels que  la forme quadratique $ax^2+bxy+cy^2$ soit primitive positive r\'eduite de discriminant~$D$ (\cite[Th.~7.2.14]{Coh93}). Dans le cas où $D=-36$, on a exactement deux repr\'esentants des classes d'\'equivalence de telles formes quadratiques donn\'es par
\begin{equation*}
x^2+9y^2\quad\textrm{et}\quad 2x^2+2xy+5y^2.
\end{equation*}
On v\'erifie alors que l'on a 
\begin{align*}
\Phi_{-36}(X)	& =(X-j(3\sqrt{-1}))\left(X-j\left(\frac{-1+3\sqrt{-1}}{2}\right)\right)\\
		& =X^2-153542016X-1790957481984=P(X)
\end{align*}
et
\begin{equation*}
j=j\left(\frac{-1+3\sqrt{-1}}{2}\right).
\end{equation*}
Cela établit l'assertion.

Par ailleurs, d'après l'expression des coefficients~$c_4$, $c_6$ et $\Delta$ ci-dessus, $E$ a réduction additive en l'idéal $\ideal_3=\sqrt{3}\integers$ et
\begin{equation*}
(v_{\ideal_3}(c_4),v_{\ideal_3}(c_6),v_{\ideal_3}(\Delta))=(4,6,9).
\end{equation*}
En particulier, on a d'après~\cite[th.~1]{Kra90}, $\lvert\Phi_{\ideal_3}\rvert=4$ ou $12$. Donc, d'après le cor.~\ref{cor:Phi_3}, la représentation $\rep$ est irréductible pour tout nombre premier $p\geq 5$. Par ailleurs, la courbe modulaire $X_0(3)$ paramétrisant les courbes elliptiques munies d'un sous-groupe stable d'ordre~$3$ est de genre~$0$ et un isomorphisme avec $\mathbf{P}^1$ est donné par la fraction rationnelle suivante~:
\begin{equation*}
j(X_0(3))(X)=\frac{(X+3)^3(X+27)}{X}.
\end{equation*}
On vérifie alors que l'on a $j(X_0(3))(243-162\sqrt{3})=j$. Cela montre que $3$ est exceptionnel. On vérifie enfin que le point de coordonnées affines
\begin{equation*}
\left\{\begin{array}{lcl}
 x & = & -2^2\cdot 7(15+8\sqrt{3}) \\
 y & = & 2^3\cdot 3\cdot 7^2(13+4\sqrt{3})
\end{array}
\right.
\end{equation*}
est un point d'ordre~$2$ de $E$. En particulier, $2$ est exceptionnel. D'où le résultat.

\medskip

\noindent{\it Remarque. }Pour montrer la finitude de l'ensemble exceptionnel, on aurait pu appliquer directement le cor.~\ref{cor:discriminant} (au lieu du cor.~\ref{cor:Phi_3}) par exemple avec le nombre premier $\ell=5$ pour lequel $B_{5}\not=0$.

\medskip

Dans les deux derniers exemples suivants, on a utilis\'e \verb"magma" (\cite{Mag97}) pour calculer certains coefficients de la fonction $L$ de~$E$.

\begin{example}
On considère $K=\Qset\left(\cos\left(\frac{2\pi}{9}\right)\right)$ le corps cubique cyclique de conducteur $9$ et la courbe $E$ d'\'equation
\begin{equation*}
y^2=x^3+2(1+\alpha)^2x+24\alpha(2+\alpha),
\end{equation*}
où $\alpha=2\cos\left(\frac{2\pi}{9}\right)$ est racine du polynôme $X^3-3X+1$. Alors, l'ensemble $\Exc[E/K]$ est vide.
\end{example}
\noindent{\it D\'emonstration. }On a $D_K=3^4$ et $3\integers=\ideal_3^3$ où $\ideal_3$ est l'idéal de $\integers$ engendré par $1+\alpha$. De plus, on vérifie que l'on a
\begin{equation*}
\left\{
\begin{array}{rl}
c_4 & =-2^5(1+\alpha)^5(1+\alpha-\alpha^2); \\
c_6 & =-2^8(1+\alpha)^{12}(1+\alpha-\alpha^2)^3; \\
\Delta & =-2^9(1+\alpha)^6\cdot 5\cdot 11 \\
\end{array}
\right.
\end{equation*}
et $1+\alpha-\alpha^2$ est une unité de $\integers$. La courbe $E$ a mauvaise réduction multiplicative en les idéaux premiers $5\integers$ et $11\integers$. En l'idéal premier $\ideal_3$, la courbe a mauvaise réduction additive avec potentiellement bonne réduction et $v_{\ideal_3}(\Delta_{\ideal_3})=v_{\ideal_3}(\Delta)=6$. Donc, d'après~\cite[th.1]{Kra90}, on a $\lvert\Phi_{\ideal_3}\rvert=2$ ou $6$. En l'idéal premier~$2\integers$, la courbe $E$ a mauvaise réduction additive avec potentiellement bonne réduction et
\begin{equation*}
(v_2(c_4),v_2(c_6),v_2(\Delta))=(5,8,9)
\end{equation*}
d'où $v_2(j)=6$ et $2v_2(c_6)=3v_2(c_4)+1$. Comme de plus
\begin{equation*}
j'=\frac{j}{2^6}=\frac{(1+\alpha)^9(1+\alpha-\alpha^2)^3}{5\cdot 11}=\frac{3^3}{5\cdot 11}\equiv 1\pmod{4},
\end{equation*}
d'après~\cite{Cal04}, on a $\lvert\Phi_2\rvert=4$. Partout ailleurs, la courbe $E$ a bonne réduction.

Vu l'\'etude pr\'ec\'edente, aucun des r\'esultats \og uniformes\fg\ du \S\ref{s:enonces_resultats} ne s'applique. Les nombres premiers $17$, $19$, $37$ et $53$ sont (totalement) décomposés dans $K$ et on vérifie que l'on a
\begin{align*}
\{t_{\ideal}\}_{\ideal\mid 17}&=\{-3,-3,3\};\qquad \{t_{\ideal}\}_{\ideal\mid 19}=\{-5,-5,5\}; \\
\{t_{\ideal}\}_{\ideal\mid 37}&=\{-7,-7,7\};\qquad \{t_{\ideal}\}_{\ideal\mid 53}=\{-3,3,3\}.
\end{align*}
Soit $p$ un nombre premier exceptionnel $\geq 5$. D'après le th.~\ref{th:critere1}, on a en particulier,
\begin{equation*}
p\in \mathfrak{S}_{17}^{(3)}\cap \mathfrak{S}_{19}^{(3)}\cap \mathfrak{S}_{37}^{(3)}=\{2,3,5\}.
\end{equation*}
Il ne reste donc plus qu'à traiter les cas où $p=2,3$ ou $5$. Or, si $\ideal_{53}$ désigne un idéal premier de $\integers$ au-dessus de $53$, le polynôme $P_{\ideal_{53}}$ est irréductible modulo~$2$ et~$5$. Par ailleurs, l'idéal $7\integers$ est premier et $t_{7}=-36$, donc le polynôme 
\begin{equation*}
P_{7}(X)=X^2+36X+7^3
\end{equation*}
est irréductible modulo~$3$. On en déduit le résultat annoncé.

\begin{example}\label{ex:biquad}
On considère $K=\Qset(\sqrt{-3},\sqrt{-7})$  et $E$ la courbe d'\'equation
\begin{equation*}
y^2=x^3+a_4x+a_6
\end{equation*}
o\`u
\begin{equation*}
\left\{\begin{array}{rcl}
 a_4 & = & \displaystyle{\frac{81}{4}}\cdot \left(69+43\sqrt{-3}+29\sqrt{-7}+17\sqrt{21}\right); \\
 & & \\
 a_6 & = & 162\cdot (207-84\sqrt{-3}-54\sqrt{-7}+46\sqrt{21}).
\end{array}
\right.
\end{equation*}
Alors, $\Exc[E/K]=\{2,3,5\}$.
\end{example}
\noindent{\it D\'emonstration. }La courbe $E$ a la propri\'et\'e particuli\`ere d'\^etre une $\Qset$-courbe, c'est-\`a-dire, d'\^etre isog\`ene \`a ses conjugu\'ees galoisiennes. Plus pr\'ecis\'ement, c'est une $\Qset$-courbe (sans multiplication complexe) de conducteur~$2\integers$ dont tous les endomorphismes sont d\'efinis sur~$K$ (cf. \cite[ex.~13]{GoGu10} et~\cite{GoLa01}). En les id\'eaux au-dessus de~$2$, elle a mauvaise r\'eduction multiplicative. En particulier, aucun des r\'esultats \og uniformes\fg\ du \S\ref{s:enonces_resultats} ne s'applique. Montrons \`a pr\'esent que le crit\`ere du th.~\ref{th:critere1} est lui aussi insuffisant pour traiter cette courbe. On doit montrer que pour tout nombre premier~$\ell\geq 3$, l'entier $B_{\ell}$ est nul, autrement dit, que $\ell^{24}$ est racine de~$P_{\ell}^*$ (lemme~\ref{lem:polynome_P_ell}). D'apr\`es les propri\'et\'es de $E$ rappel\'ees ci-dessus, on a $P_{\ideal}=P_{\ideal'}$ pour tout couple $(\ideal,\ideal')$ d'id\'eaux premiers divisant~$\ell$. Or $D_K=3^2\cdot 7^2$, donc si $\ell\not=3,7$, $\ell\integers$ se d\'ecompose en un produit de $2$ ou $4$ id\'eaux premiers. On a alors respectivement
\begin{equation*}
P_{\ell}^*=\left(P_{\ideal}^{(12)}\right)^{*2}\quad\textrm{et}\quad P_{\ell}^*=\left(P_{\ideal}^{(12)}\right)^{*4}.
\end{equation*}
Or, dans le premier cas, les racines complexes $\alpha$ et $\beta$ de $P_{\ideal}$ satisfont $\alpha\beta=\ell^2$ et dans le second, $\alpha\beta=\ell$. On en d\'eduit le r\'esultat voulu dans ce cas. Par ailleurs, on v\'erifie que l'on a $P_3^*(X)=(X-3^{24})^2$ et
\begin{equation*}
P_7^*(X)=(X-7^{24})^2\cdot (X^2-371728108602950083202X+7^{48}),
\end{equation*}
d'o\`u la nullit\'e de $B_{\ell}$ pour tout~$\ell$ de bonne r\'eduction. Pour cette courbe, on a donc recours au crit\`ere du th.~\ref{th:critere2}. Le nombre de classes $h$ de~$K$ est~$1$. On consid\`ere l'id\'eal premier~$\ideal_5$ au-dessus de~$5$ engendr\'e par une racine~$\gamma_{\ideal_5}$ du polyn\^ome $\mathfrak{m}_{\gamma_{\ideal_5}}(X)=X^4+17X^2+25$. On a alors,
\begin{equation*}
P_{\ideal_5}(X)=X^2+4X+25\quad\textrm{d'o\`u}\quad P_{\ideal_5}^{(12)}(X)=X^2 -2\cdot 47\cdot 1163039X + 5^{24}
\end{equation*}
et
\begin{equation*}
\mathfrak{m}_{\gamma_{\ideal_5}}^{(12)}(X)=\left(X^2-2\cdot 73\cdot 19441X+5^{12}\right)^2
\end{equation*}
puis,
\begin{equation*}
\left(\mathfrak{m}_{\gamma_{\ideal_5}}^{(12)}\right)^{* 2}(X)=\left(X-5^{12}\right)^8\cdot \left(X^2-2\cdot 79\cdot 127\cdot 337 \cdot 1191313X+5^{24}\right)^4.
\end{equation*}
On en d\'eduit que l'on a 
\begin{multline*}
R_{\ideal_5}=2^{126}\cdot 3^{100}\cdot 5^{225}\cdot 7^{21}\cdot 11\cdot 13^{8}\cdot 19\cdot 37^{11}\cdot 41^{8}\cdot 59^{16}\cdot 103\cdot 109^{8}\cdot 149^{8} \\
\cdot 193\cdot 373^{2}\cdot 2137\cdot 4201^{2}\cdot 7753^{2}\cdot 24061^{2}.
\end{multline*}
On recommence ensuite ces m\^emes calculs avec l'id\'eal premier~$\ideal_7$ au-dessus de~$7$ engendr\'e par une racine~$\gamma_{\ideal_7}$ du polyn\^ome $\mathfrak{m}_{\gamma_{\ideal_7}}(X)=X^4+4X^3+11X^2+14X+7$. On a alors $P_{\ideal_7}(X)=X^2+2X+7$ et
\begin{multline*}
R_{\ideal_7}=2^{105}\cdot 3^{59}\cdot 5^{26}\cdot 7^{116}\cdot 11^{2}\cdot 13^{2}\cdot 17^8\cdot 23^{8}\cdot 31\cdot 79\cdot 137^2\cdot 191^{4}\cdot 193 \\
\cdot 463\cdot 487^{2}\cdot 673\cdot 1033^{2}\cdot 1471\cdot 2953\cdot 3697.
\end{multline*}
Apr\`es ces deux it\'erations du th.~\ref{th:critere2}, on a donc montr\'e l'inclusion
\begin{equation*}
\Exc[E/K]\subset\{2,3,5,7,11,13,193\}.
\end{equation*}
Notons respectivement, $\ideal_3$ et $\ideal_{17}$ un id\'eal premier au-dessus de~$3$ et de~$17$. Alors, le polyn\^ome $P_{\ideal_3}(X)=X^2+9$ est irr\'eductible modulo~$7$ et $11$ et le polyn\^ome $P_{\ideal_{17}}(X)=X^2+10X+289$ est irr\'eductible modulo~$193$. De m\^eme, le polyn\^ome $P_{\ideal_{5}}$ ci-dessus est irr\'eductible modulo~$13$. Enfin, $2,3$ et $5$ sont exceptionnels car ce sont les degr\'es des isog\'enies de~$E$ vers ses trois conjugu\'ees galoisiennes (loc.~cit.). D'o\`u le r\'esultat.

\providecommand{\bysame}{\leavevmode\hbox to3em{\hrulefill}\thinspace}

\end{document}